\documentclass[11pt, dvips, twoside]{amsart}
\usepackage{amssymb,latexsym,bbm,graphicx,epsfig,epic,eepic,psfrag,oldgerm}

\def\aa{\alpha}
\def\bb{\beta}
\def\gg{\gamma}
\def\dd{\delta}
\def\ee{\epsilon}

\def\ff{\varphi}

\def\oo{\omega}
\def\ss{\sigma}

\def\SS{\Sigma}

\def\ce{{\mathcal E}}

\def\cj{{\mathcal J}}
\def\cl{{\mathcal L}}

\def\ra{\rightarrow}
\def\Ra{\Rightarrow}

\def\ha{\hookrightarrow}

\def\sq{\Box}
\def\tr{\triangle}

\def\proof{\noindent {\it Proof. \;}}

\def\ni{\noindent}

\def\s{\smallskip}
\def\m{\medskip}
\def\b{\bigskip}

\newcommand{\proofend}{\hspace*{\fill} $\Box$\\}
\newcommand{\diam}{\hspace*{\fill} $\Diamond$}

\theoremstyle{plain}
\newtheorem{theorem}{Theorem}[section]
\newtheorem{lemma}[theorem]{Lemma}
\newtheorem{proposition}[theorem]{Proposition}
\newtheorem{corollary}[theorem]{Corollary}
\newtheorem{remark}[theorem]{Remark}

\newtheorem{definition}[theorem]{Definition}

\newtheorem{question}[theorem]{Question}

\def\Vol{\operatorname {vol}\:\!}

\def\Diffc0{\operatorname{Diff^c_0}}

\def\Sympc0{\operatorname{Symp^c_0}}

\def\Int{\operatorname{Int}}

\def\id{\operatorname{id}}
\def\pr{\operatorname{pr}}
\def\im{\operatorname{Im}\:\!}
\def\Vol{\operatorname{Vol}\;\!}

\def\Int{\operatorname{Int}\:\!}
\def\conv{\operatorname{conv}}

\def\sym{\operatorname{Sp}(n;\RR)}

\def\diag{\operatorname{diag}}

\def\CC{\mathbbm{C}}

\def\NN{\mathbbm{N}}
\def\PP{\mathbbm{P}}
\def\QQ{\mathbbm{Q}}
\def\RR{\mathbbm{R}}

\def\ZZ{\mathbbm{Z}}

\begin{document}

\title{Packing symplectic manifolds by hand}

\author{Felix Schlenk}
\address{(F.\ Schlenk) Mathematisches Institut,
Universit\"at Leipzig, 04109 Leipzig, Germany}
\email{schlenk@math.uni-leipzig.de}

\date{\today}
\thanks{This work is supported by the Deutsche Forschungsgemeinschaft}

\begin{abstract}
We construct explicit maximal symplectic packings of 
minimal rational and ruled symplectic $4$-manifolds by few balls
in a very simple way.
\end{abstract}

\maketitle

\markboth{
{\rm Packing symplectic manifolds by hand}
}{{}} 


\setcounter{tocdepth}{2}

\setcounter{secnumdepth}{4}

\section{Introduction}

\ni
Consider a connected $2n$-dimensional symplectic
manifold $(M, \oo)$ of finite volume $\Vol (M, \oo) =
\frac{1}{n!} \int_M \oo^n$, and let $B^{2n}(a)$ be the open ball
of radius $\sqrt{a / \pi}$ in standard symplectic space $\left(
\RR^{2n}, \oo_0 \right)$.
The {\it $k$'th symplectic packing number} $p_k (M, \oo) \in \;]0,1]$
is defined as
\[
p_k (M, \oo) \,=\, \sup_a \frac{k \Vol \left( B^{2n}(a),\oo_0 \right)}{\Vol
(M, \oo) }  
\]
where the supremum is taken over all those $a$ for which the
disjoint union $\coprod_{i=1}^k B^{2n}(a)$ of $k$ equal balls
symplectically embeds into $(M, \oo)$.
%
%
If $p_k(M, \oo) <1$, one says that there is a {\it packing
obstruction},
and if $p_k(M, \oo) =1$, one says that $(M,\oo)$ admits a 
{\it full packing}\, by $k$ balls.
The first examples of packing obstructions were found by Gromov,
\cite{G1},
and many further packing obstructions and also some exact values of $p_k$
were obtained by Mc\:\!Duff and Polterovich in \cite{MP}. 
Finally, Biran showed in \cite{B1,B2} that 
\begin{equation}  \label{P:finite}
P(M, \oo) \,:=\, 
\inf \left\{ k_0 \in \NN \mid p_k(M,\oo) =1 \quad \mbox{for all
}\, k \ge k_0 \right\} \,<\, \infty
\end{equation}
for an interesting class of closed symplectic $4$-manifolds
containing sphere bundles over a surface 
and for all closed symplectic $4$-manifolds with $[\oo] \in
H^2(M;\QQ)$.

Besides sporadic results on the first packing number $p_1$
and besides the determination of $p_2 \left( E(a_1, \dots, a_n) \right)$
in \cite{MMT},
all known computations of packing numbers are contained in 
\cite{MP,B1,B2}.
We refer to Biran's excellent survey \cite{B-ECM} for the methods used, 
and only mention that in \cite{MP,B1,B2} the problem of
symplectically embedding $k$ equal balls into $(M,\oo)$ is first
reformulated as the problem of deforming a symplectic form 
on the $k$-fold blow-up of $(M,\oo)$ along a
certain family of cohomology classes,
and that this problem is then solved using tools from classical
algebraic geometry, Seiberg--Witten--Taubes theory, and
Donaldson's symplectic submanifold theorem, respectively.
As a consequence, the symplectic packings found are not explicit.
For some of the symplectic manifolds considered in
\cite{MP,B1,B2} and some values of $k$, 
explicit maximal symplectic packings were constructed by Karshon
\cite{Ka},
Traynor \cite{T}, Kruglikov \cite{Kr}, and Maley, Mastrangeli 
and Traynor, \cite{MMT}.
In this article we construct all known and also some new  
explicit maximal packings of symplectic $4$-manifolds in a very simple
way.
To be more precise, we construct maximal packings different from
those in \cite{Ka,T,Kr,MMT} of the $4$-ball and of $\CC \PP^2$ by
$k \le 6$ balls and by $l^2$ balls for each $l \in \NN$, 
of the product of two surfaces of
equal area by $2l^2$ balls, and of the ellipsoids $E(\pi,
k\pi)$ and $E(\pi,a)$ by $k$ and $2$ balls, respectively.
In addition, we construct maximal packings of $S^2 \times S^2$
by $k \le 6$ balls for all symplectic structures and by $7$
balls for some symplectic structures,
as well as maximal packings of the non-trivial bundle $S^2 \ltimes S^2$
by $k \le 5$ balls for all symplectic structures and by $6$
balls for some symplectic structures.
In the range of $k$ for which these constructions fail to give
maximal packings, they give a feel that the balls in the packings from
\cite{MP,B1,B2} must be ``wild''.
We shall also construct an explicit full packing of the $2n$-ball by $l^n$
equal balls for each $l \in \NN$ in a most simple way.

In the next section we give several motivations for the
symplectic packing problem.
In Section~\ref{s:numbers} we collect the packing numbers of
interest to us, and in Section~\ref{s:construction} we construct
our maximal packings of symplectic $4$-manifolds. In the last
section we overview what is known in dimensions $\ge 6$ and
construct full packings of the $2n$-ball by $l^n$ balls.

Balls will always be endowed with the standard symplectic form 
$\oo_0 = \sum_{i=1}^n dx_i \wedge dy_i$. 
Since the packing numbers $p_k \left( B^{2n} (a), \oo_0 \right)$
do not depend on $a$, we shall usually pack the unit ball
$B^{2n} := \left( B^{2n}(\pi), \oo_0 \right)$.

\section{Motivations for the symplectic packing problem}

\ni
{\bf 1. Higher Gromov widths}

\s
\ni
The Gromov width 
\[
w_G (M,\oo) \,:=\, 
  \sup \{\:\! a \mid B^{2n}(a) \text{ symplectically embeds into
} (M, \oo) \:\!\} 
\]
of a symplectic manifold $(M, \oo)$ measures the size of a
largest Darboux chart of $(M,\oo)$. 
It is the smallest normalized symplectic capacity as defined
in \cite[p.\ 51]{HZ}, and we refer to
\cite{B1,B2,BC,FGS,G1,J,KT,LM1,LM2,Lu,M0,M-variants,MP,MSl,buch,Sib} 
for results on the Gromov width and to
Section~\ref{s9:3} below for explicit symplectic
embeddings realizing $w_G(M,\oo)$ or estimating it from below. 
If $(M,\oo)$ has finite volume,
the first packing number $p_1(M,\oo)$ is equivalent to the Gromov 
width,
\[
p_1 (M, \omega) \Vol (M, \oo) \,=\, \frac{1}{n!} \bigl( w_G (M,\oo) \bigr)^n .
\] 
Similarly, the higher packing numbers $p_k(M,\oo)$, $k \ge 2$,
are equivalent to the {\it higher Gromov widths}
\[
w_G^k (M,\oo) \,:=\, 
  \sup \left\{ a \,\Bigg|\, \coprod_{i=1}^k B^{2n}(a) \text{
symplectically embeds into } (M, \oo) \right\} ,
\]
which form a distinguished sequence of embedding capacities as
considered in \cite{CHLS}.

\m
\ni
{\bf 2. ``Superrecurrence for symplectomorphisms'' via packing obstructions?}

\s
\ni
In view of the Poincar\'e recurrence theorem, volume preserving mappings
have strong recurrence properties. 
The solution of the Arnold conjecture for the torus by Conley and 
Zehnder \cite{CZ} in 1983 
demonstrated that Hamiltonian
symplectomorphisms have yet stronger recurrence properties. 
As was pointed out to me by Polterovich, the original motivation for Gromov 
to study the packing numbers $p_k$ was his search for
recurrence properties of arbitrary symplectomorphisms which are stronger
than those of volume preserving mappings.

We explain the relation between ``superrecurrence for
symplectomorphisms'' and symplectic packing obstructions by means
of an example.
Let $B$ and $B'$ be the open balls in $\RR^{2n}$ centred at the
origin of volumes $2^n - \frac 12$ and $1$, respectively.
For every compactly supported volume preserving diffeomorphism
$\ff$ of $B$ set 
\[
R(\ff) \,=\, \min \left\{ m \in \NN \mid \ff^m \left( B' \right)
\cap B' = \emptyset \right\} .
\]
Of course, $R(\ff) \le 2^n -2$, and using Moser's deformation
argument, for which we refer to \cite[p.\ 11]{HZ}, it is easy to
construct a $\ff$ with $R(\ff) = 2^n-2$. 
The packing obstruction $p_2(B) = \frac{1}{2^{n-1}}$ proved by
Gromov in \cite{G1} shows, however, that $R(\ff) =1$ if $\ff$ is symplectic.

This motivation for symplectic packings lost some of its appeal
by the work of Mc\;\!Duff--Polterovich and Biran.
Indeed, in dynamics one usually asks for recurrence into {\it
small}\, neighbourhoods of a point. To establish recurrence of
small balls we would need packing obstructions for {\it
large}\:\! $k$. 
In view of \cite[Remark~1.5.G]{MP}, these obstructions asymptotically
always vanish, and in view of \eqref{P:finite}, they completely
vanish for many symplectic $4$-manifolds.

\m
\ni
{\bf 3. Between Euclidean and volume preserving}

\s
\ni
{\it Volume preserving packings.}
Consider a connected $n$-dimensional manifold $M$ endowed with a
volume form $\Omega$ such that the volume $\Vol (M, \Omega) =
\int_M \Omega$ is finite,
and denote the Lebesgue measure of an open subset $U$ of $\RR^n$
by $\left| U \right|$.
We write $B^n(A)$ for the open ball of radius $\sqrt{A/\pi}$ in $\RR^n$.
For $k \in \NN$ we set
\[
v_k (M, \Omega) \,=\, \sup \left\{ \frac{k \left| B^n (A)
\right|}{\Vol (M,\Omega)} \right\}
\]
where the supremum is taken over all $A$ for which there exists
a volume preserving embedding $\coprod_{i=1}^k B^{n}(A) \hookrightarrow (M,\Omega)$.
Moser's deformation method readily implies that
$v_k (M, \Omega) =1$ for all $k \in \NN$.
The main result of \cite{Sch} shows more:
For any partition $M = \coprod_{i=1}^k M_i$ of $M$ into subsets
$M_i$
such that $\Int M_i$ is connected and $\Vol \left( \Int M_i,
\Omega \right) = \frac 1k \Vol (M, \Omega)$ for all $i$ there
exists a volume preserving embedding 
$\coprod_{i=1}^k B^n(A) \ha \coprod_{i=1}^k \Int M_i$ with
$\left| B^n(A) \right| = \frac 1k \Vol (M, \Omega)$.
If the volume form $\Omega$ comes from a symplectic form $\oo$, 
the sequence $\left(1 - p_k (M,\omega) \right)_{k \in \NN}$ is a
measure for how far the symplectic geometry of $(M,\oo)$ is from
the volume geometry of $(M, \Omega)$.

\m
\ni
{\it Euclidean packings.}
Given a bounded domain $U$ in $\RR^n$, define its 
{\it $k$'th Euclidean packing number}\, as 
\[
\Delta_k (U) \,=\, \sup \left\{ \frac{k \left| B^n(a)
\right|}{\left| U \right|}\right\}
\]
where the supremum is taken over all $a$ for which $k$ disjoint
translates of $B^n(a)$ fit into $U$. Then
\[
\Delta_k(U) \le p_k(U) \le v_k(U) =1 \quad \text{ for all }\, k
\in \NN ,
\]
and it is interesting to understand ``on which side'' $p_k(U)$ lies.
To fix the ideas, we assume that $U$ is the unit ball 
$B^n := B^n(\pi)$ in $\RR^n$.
The precise values of $\Delta_k \left( B^n \right)$ are known
only for small $k$:
If $1 \le k \le n+1$, the smallest ball containing $k$ balls of
radius $1$ has radius $1+\sqrt{2-2/k}$, and the centres of
the balls are arranged as vertices of a regular
$(k-1)$-dimensional simplex inscribed in the ball and concentric
with it.
Moreover, if $n+2 \le k \le 2n$, the smallest ball $B$
containing $k$ balls of radius $1$ has radius 
$1+ \sqrt{2}$, and the packing configuration of $2n$ balls in
$B$ is unique up to isometry, the centres being the midpoints of
the faces of an $n$-dimensional Euclidean cube whose edges have
length $2 \sqrt{2}$.
In particular,
\begin{equation}  \label{Deltak}
\Delta_k \left( B^n \right) \,=\, 
\left\{
  \begin{array}{ccl}   
         \frac{k}{\left( 1 + \sqrt{2 - \frac 2 k}\right)^n}
                   & \text{ if }\, & 1 \le k \le n+1 ,  \\ [.8em]
         \frac{k}{\left( 1 + \sqrt{2} \right)^n}
                   & \text{ if }\, & n+2 \le k \le 2n .
  \end{array}   
\right.                   
\end{equation}
While for $1 \le k \le n+1$ these numbers were known to Rankin
in 1955, for $n+2 \le k \le 2n$ they  were obtained only
recently by W.\ Kuperberg, \cite{Ku}.
An obvious upper bound for $\Delta_k \left( B^n \right)$ is
\begin{equation}  \label{est9:obvious}
\Delta_k \left( B^n \right) \,\le\, \frac{k}{2^n}
\quad\; \text{ for all }\, k \ge 2 . 
\end{equation}
Given a bounded domain $U$ in $\RR^n$, let $\conv (U)$ be
the convex hull of $U$.
For each $k \ge 1$ we set
\begin{equation}  \label{def:conv}
\conv_k \left( B^n \right) \,=\, 
\sup \frac{ k \left| B^n \right|}{\left| \conv (U) \right|}
\end{equation}
where the supremum is taken over all configurations $U$ of $k$ disjoint
translates of $B^n$ in $\RR^n$.
Since $B^n$ is convex,
$\Delta_k \left( B^n \right) \le \conv_k \left( B^n \right)$ for all $k
\in \NN$.
Let $S_k^n = \conv (U)$ be the sausage
obtained by choosing
\begin{equation}  \label{eq9:U=}
U \,=\, \coprod_{i=0}^{k-1} \left( B^n + i \mathbf{u} \right)
\end{equation}
where $\mathbf{u}$ is a unit vector in $\RR^n$.
With $\kappa_n := \left| B^n \right|$ we then have 
$\left| S_k^n \right| = \kappa_n +
2(k-1) \kappa_{n-1}$.
The sausage conjecture of L.\ Fejes T\'oth from 1975 states that equality in
\eqref{def:conv} is attained exactly for $U$ as in \eqref{eq9:U=}, and
this conjecture was proved by Betke and Henk, \cite{BH}, for $n \ge 42$.
Therefore, 
\begin{equation}  \label{ine9:sausage}
\Delta_k \left( B^n \right) \,\le\, \conv_k \left( B^n \right) \,=\,
\frac{k \:\! \kappa_n}{\kappa_n + 2 (k-1) \kappa_{n-1}}
\,<\, \frac{k}{k-1} \sqrt{\frac \pi 2} \sqrt{\frac{1}{n+1}}
\,\quad \text{ if }\; n \ge 42 .
\end{equation}
For arbitrary $n$, an older result of Gritzmann, \cite{Gr}, 
states that
\[
\Delta_k \left( B^n \right) \,\le\, \conv_k \left( B^n \right) \,<\,
(2+\sqrt{3}) \sqrt{\frac{\pi}{2}} \sqrt{\frac{1}{n}} .
\]
In order to get an idea of the values  
$\Delta_k \left( B^n \right)$ for large $k$ we notice that the
limit
\[
\Delta^n \,:=\, \lim_{k \ra \infty} \Delta_k \left( B^n \right)
\]
exists and is equal to the highest density of a packing of 
$\RR^n$, see Section~2.1 of Chapter~3.3 in \cite{Handbook}. 
The highest density of a packing of $\RR^2$ is 
\[
\Delta^2 \,=\, \frac{\pi}{\sqrt{12}} \,=\, 0.9069\dots
\]
as in the familiar hexagonal lattice packing in which each disk
touches 6 others (Thue, 1910).
The highest density of a packing of $\RR^3$ is 
\[
\Delta^3 \,=\, \frac{\pi}{\sqrt{18}} \,=\, 0.74048\dots
\]
as in the face centred cubic lattice packing which is 
usually found in fruit stands and in which each ball touches 12
other balls.
This was conjectured by Keppler in 1611, and Gauss proved in
1831 that no lattice packing has a higher density.
The Keppler conjecture 
was settled only recently by Hales, see
\cite{Ha} and the references therein.
For $4 \le n \le 36$, the currently best upper bound for $\Delta^n$
was given recently by Cohn and Elkies in \cite{CE}. 
E.g., 
\[
\frac{\pi^2}{16} = 0.61685 \le \Delta^4 \le 0.647742 .
\]
Here, the lower bound is the density of the packing associated
with the ``checkerboard lattice'' consisting of all vectors
$(a,b,c,d) \in \ZZ^4$ with $a+b+c+d \in 2 \ZZ$, and it is known
that this is the highest possible density for a $4$-dimensional 
lattice packing.
A result of Blichfeldt from 1929 states that
\begin{equation}  \label{est:Blichfeldt}
\Delta^n \,\le\, (n+2) 2^{-(n+2)/2} , 
\end{equation}
and the best known lower and upper bounds for $\Delta^n$ of
asymptotic nature are
\[
c\:\! n 2^{-n} \,\le\, \Delta^n \,\le\, 2^{- \left(0.599+o(1)\right) n}
\,\quad \text{ as }\, n \ra \infty 
\]
for any constant $c < \log 2$, see Section~2 of Chapter~3.3 in 
\cite{Handbook}.

We refer to \cite{CS}, to Sections~3.3 and 3.4 of \cite{Handbook}, and
to \cite{Zong}
for more information on Euclidean packings, its long history and
its many relations and applications to 
other branches of mathematics (such as discrete geometry, 
group theory, number theory and crystallography) and to problems in physics,
chemistry, engineering and computer science.

The symplectic packing numbers $p_k \left( B^4 \right)$ are listed in 
Table~\ref{ta:ball} below. 
For $n \ge 3$, the results known about $p_k \left( B^{2n} \right)$ are
\begin{eqnarray}     
   p_k \left( B^{2n} \right) &=& \frac{k}{2^n}  \quad\, \text{ for }\,
                        2 \le k \le 2^n,  \label{eq9:pkBkn} \\ [.3em]
   p_{l^n} \left( B^{2n} \right) &=& 1  \quad\, \text{ for all }\,  
                         l \in \NN, \label{eq9:pkBln} 
\end{eqnarray}
see \cite[Corollary~1.5.C and 1.6.B]{MP} and Section~\ref{s9:higher}.1 below. 
The identities~\eqref{eq9:pkBln} show that
\begin{equation}  \label{ball9:asym}
\lim_{k \ra \infty} p_k \left( B^{2n} \right) =1 \quad\, \text{
for all }\; n.
\end{equation}

Of course, $\Delta_k \left( B^2 \right) < p_k \left( B^2 \right)
= v_k \left( B^2 \right) =1$ for all $k \ge 2$. 
Comparing \eqref{Deltak} or \eqref{est9:obvious} for $n=4$ 
with the values $p_k \left( B^4 \right)$ listed in Table~\ref{ta:ball} 
we see that 
\[
\Delta_k \left( B^4 \right) \,<\, p_k \left( B^4 \right)
\quad\, \text{ for all }\; k \ge 2 .
\]
Moreover, \eqref{est9:obvious} and \eqref{eq9:pkBkn} show that
\[
\Delta_k \left( B^{2n} \right) \,\le\, \frac{1}{2^n}\, p_k \left( B^{2n}
\right) 
\quad\, \text{ for }\; 2 \le k \le 2^n \;\text{ and all }\; n \in \NN .
\]
Inequality~\eqref{ine9:sausage} and \eqref{eq9:pkBln} yield an explicit
$k(2n)$ such that
\[
\Delta_k \left( B^{2n} \right) \,<\, p_k \left( B^{2n} \right)
\quad\, \text{ for all }\; k \ge k(2n) \;\text{ and }\; 2n \ge 42 .
\]
It is conceivable that $\Delta_k \left( B^{2n} \right) < p_k \left(
B^{2n} \right)$ for all $k \ge 2$ and $n \in \NN$, but we do not know
the answer to
\begin{question}
Is it true that $\Delta_{28} \left( B^6 \right) < p_{28} \left(B^6 \right)$?
\end{question}
Finally, comparing \eqref{ball9:asym} with \eqref{est:Blichfeldt} we see that
$p_k \left( B^{2n} \right)$ is much larger than $\Delta_k \left( B^{2n} \right)$
for sufficiently large $k$ and large $n$.

\m
\ni  
{\bf 4. Relations to algebraic geometry}

\s
\ni
A symplectic packing of $(M, \oo)$ by $k$ equal balls
corresponds to a symplectic blow up of $(M,\oo)$ at $k$ points
with equal weights. 
Via this correspondence, the symplectic packing problem is
intimately related to old problems in algebraic geometry:
The symplectic packing problem for the complex
projective plane $\CC \PP^2$ (completely solved by Biran in \cite{B1})
is related to an old (and still open) conjecture of Nagata on the
minimal degree of an irreducible algebraic curve in $\CC \PP^2$
passing through $N \ge 9$ points with given multiplicities, see
\cite{B-ample,B-ECM,MP,Xu} for details.
Moreover, the symplectic packing problem is closely related to
the problem of computing Seshadri constants of ample line
bundles, which are a measure of their local positivity, see
\cite{B-ample,B-ECM,BC,La}.


\section{The packing numbers of the $4$-ball, of $\CC \PP^2$
and of ruled symplectic $4$-manifolds} 
\label{s:numbers}

\ni
In this section we review the known packing numbers of interest to us and also
compute $p_k$ for the nontrivial sphere bundles over Riemann surfaces
for $k \le 7$. 

\subsection{The packing numbers of the $4$-ball and of $\CC \PP^2$}  \label{11}

Let $\oo_{SF}$ be the unique $\mbox{U}(3)$-invariant K\"ahler form on 
$\CC \PP^2$ whose integral over $\CC \PP^1$ equals $1$. 
According to a  result of Taubes, \cite{Ta}, 
every symplectic form on $\CC \PP^2$ is diffeomorphic to $a\,
\oo_{SF}$ for some $a \neq 0$.
In view of the symplectomorphism
\begin{equation}  \label{emb:BCP2}
\left( B^4(\pi), \oo_0 \right) \ra \left( \CC \PP^2 \setminus
\CC \PP^1, \pi \,\oo_{SF} \right), \quad\;
{\boldmath{z}} = (z_1, z_2) \mapsto \left[ z_1 : z_2 : 
                \sqrt{1- | {\boldmath{z}} |^2} \, \right]
\end{equation}
further discussed in \cite[Example~7.14]{MS}
we have $p_k \left( B^4 \right) \le p_k \left( \CC \PP^2 \right)$ for all $k$. 
It is shown in \cite[Remark~2.1.E]{MP} that in fact
\begin{equation}  \label{pk=}
p_k \left( B^4 \right) \,=\, p_k \left( \CC \PP^2 \right)
\quad\; \text{for all }\, k .
\end{equation}
A complete list of these packing numbers was obtained in
\cite{B1} (see Table~\ref{ta:ball}).

\begin{table}[h]  
 \begin{center}
 \renewcommand{\arraystretch}{1.5}
  \begin{tabular}{cccccccccc} \hline
  $k$    & 1 &      2      &     3         & 4 &      5        &     6
  & 7 & 8 & $\ge 9$      \\ \hline
  $p_k$  & 1 & $\frac{1}{2}$ & $\frac{3}{4}$ & 1 & $\frac{20}{25}$ &
  $\frac{24}{25}$ &
  $\frac{63}{64}$ & $\frac{288}{289}$ & 1  \\ \hline
  \end{tabular}
 \end{center}
 \caption{$p_k \left( B^4 \right) = p_k \left( \CC \PP^2 \right)$} \label{ta:ball}
\end{table}

\ni
Explicit maximal packings were found by Karshon \cite{Ka}
for $k \le 3$ and by Traynor \cite{T} 
for $k \le 6$ and $k = l^2 \; (l \in \NN)$. 
We will give even simpler maximal packings for these values of $k$
in \ref{21}.

\subsection{The packing numbers of ruled symplectic $4$-manifolds}  \label{12}

Denote by $\Sigma_g$ the closed orientable surface of genus $g$. There are
exactly two orientable $S^2$-bundles with base $\Sigma_g$, namely the
trivial bundle $\pi \colon \Sigma_g \times S^2 \ra \Sigma_g$ and the
nontrivial bundle $\pi \colon \Sigma_g \ltimes S^2 \ra
\Sigma_g$, see
\cite[Lemma 6.9]{MS}. Such a manifold $M$ is called a 
{\it ruled surface}.
A symplectic form $\oo$ on a ruled surface is called {\it compatible}\,
with the given ruling $\pi$ 
if it restricts on each fibre to a
symplectic form. Such a symplectic manifold is then called a {\it ruled
symplectic $4$-manifold}. 
It is known that every symplectic structure on a
ruled surface is diffeomorphic to a form compatible with the given
ruling $\pi$ via a diffeomorphism which acts trivially on homology, and
that two cohomologous symplectic forms compatible with the same ruling
are isotopic \cite{LM4}.
A symplectic form $\oo$ on a ruled surface $M$ is thus determined up to
diffeomorphism by the class $[\oo] \in H^2(M; \RR)$.
In order to describe the set of cohomology classes realized by
(compatible) forms on $M$ 
we fix an orientation of $\Sigma_g$ and an orientation of the fibres of 
the given ruled surface $M$.
These orientations determine an orientation of $M$ in a natural way, see
below.
We say that a compatible symplectic form $\oo$ is {\it
admissible}\, if 
its restriction to each fibre induces the given orientation and if $\oo$
induces the natural orientation on $M$.
Notice that every symplectic form on $M$ is diffeomorphic to an
admissible form for a suitable choice of orientations of $\Sigma_g$ and
the fibres.

Consider first the trivial bundle $\Sigma_g \times S^2$, 
and let $\{ B = [\Sigma_g \times pt], F = [pt \times S^2]
\}$ be a basis of $H^2(M;\ZZ)$.
Here and henceforth we identify homology and cohomology
via Poincar\'e duality. 
The natural orientation of $\Sigma_g \times S^2$ is such that $B \cdot F
=1$. 
A cohomology class $C = bB + aF$ can
be represented by an admissible form if and only if 
$C \cdot F >0$ and $C \cdot C >0$, i.e.,
\[
a >0 \quad \text{and} \quad b >0 ,
\]
standard representatives being split forms.
We write $\Sigma_g(a) \times S^2(b)$ for this ruled
symplectic $4$-manifold. 

In case of the nontrivial bundle $\Sigma_g \ltimes S^2$ a basis of
$H^2(\Sigma_g \ltimes S^2; \ZZ)$ is given by $\{ A, F \}$, where $A$ is
the class of a section with selfintersection number $-1$ and $F$ is the
fibre class. 
The homology classes of sections of $\Sigma_g \ltimes S^2$ of
self-intersection number $k$ are $A_k = A + \frac{k+1}{2} F$ with $k$
odd.
The natural orientation of $\Sigma_g \ltimes S^2$ is such that $A_k
\cdot F = A \cdot F = 1$ for all $k$.
Set $B= A+ F/2$. Then $\{ B, F\}$ is a basis of
$H^2(\Sigma_g \ltimes S^2; \RR)$ with $B \cdot B = F \cdot F = 0$ and $B
\cdot F =1$. 
As for the trivial bundle, the necessary condition for a cohomology
class $bB + aF$ to be representable by an admissible form is $a>0$ and
$b>0$. 
It turns out that this condition is sufficient only if $g \ge 1$:
A cohomology class $bB + aF$ can be represented by an admissible form 
if and only if 
\begin{equation*}
      \begin{array}{lll}   
              a > b/2 >0   & \mbox{ if } & g = 0,  \\[.2em]
              \text{$a>0$ and $b>0$}  & \mbox{ if } & g \ge 1,
       \end{array}   
\end{equation*}
see \cite[Theorem 6.11]{MS}.
We write $(\Sigma_g \ltimes S^2, \oo_{ab})$ for this ruled symplectic
$4$-manifold.
A ``standard K\"ahler form'' in the class $[\oo_{ab}]$ is  explicitly
constructed in \cite[Section~3]{M1} 
and \cite[Exercise~6.14]{MS}.
When constructing our explicit symplectic packings, it will always be
clear which symplectic form in $[\oo_{ab}]$ is chosen.

We begin with the trivial sphere bundle over the sphere.
\begin{proposition}  \label{p:C21}
Assume that $a \ge b$. Abbreviate $p_k = p_k(S^2(a) \times
S^2(b))$, and denote by $\lceil x \rceil$ the minimal integer which is greater
than or equal to $x$. Then
\[
p_k = \frac{k}{2}\frac{b}{a} \quad \mbox{ if } \, \left\lceil \frac{k}{2}
\right\rceil \frac{b}{a} \le 1.
\]
Moreover,
\begin{gather*}
p_1 = \frac{b}{2a}, \qquad 
p_2 = \frac{b}{a},  \qquad 
p_3 = \frac{3}{2ab} \left\{ b, \frac{a+b}{3} \right\}^2   \;
      \mbox{ on } \left] 0, \frac{1}{2}, 1 \right],                      \\
p_4 = \frac{4}{3} p_3,  \qquad 
p_5 = \frac{5}{2ab} \left\{ b, \frac{a+2b}{5} \right\}^2  \;
      \mbox{ on } \left] 0, \frac{1}{3}, 1 \right],                      \\
p_6 = \frac{3}{ab} \left\{ b, \frac{a+2b}{5},
      \frac{2a+2b}{7} \right\}^2 \;
      \mbox{ on } \left] 0, \frac{1}{3}, \frac{3}{4}, 1 \right],         \\
p_7 = \frac{7}{2ab} \left\{ b, \frac{a+3b}{7}, \frac{3a+4b}{13},
      \frac{4a+4b}{15} \right\}^2 \; 
      \mbox{ on } \left] 0, \frac{1}{4}, \frac{8}{11}, \frac{7}{8}, 1 \right] .     
\end{gather*}
In particular, for $k \le 7$ we have $p_k(S^2(a) \times S^2(b)) =1$
exactly for 
\mbox{$(k=2, \,\frac{b}{a}=1)$}, 
$(k=4, \,\frac{b}{a}=\frac{1}{2})$,
$(k=6, \,\frac{b}{a}=\frac{1}{3})$, 
$(k=6, \,\frac{b}{a}=\frac{3}{4})$ and
\mbox{$(k=7, \,\frac{b}{a}=\frac{7}{8})$}.
\end{proposition}
We explain our notation by an example: $p_3 = \frac{3}{2ab} b^2$ if
$0 < \frac{b}{a} \le \frac{1}{2}$ 
and $p_3 = \frac{3}{2ab} \left( \frac{a + b}{3} \right)^2$ if
$\frac{1}{2} \le \frac{b}{a} \le 1$.
\\
\\
In \ref{221} we will construct explicit maximal packings of
$S^2(a) \times S^2(b)$ for 
all $k$ with $\left\lceil \frac{k}{2} \right\rceil \frac{b}{a} \le 1$,
for $k \le 6$ and $0 < b \le a$ arbitrary, and for $k=7$ and $0 <
\frac{b}{a} \le \frac{3}{5}$,  
as well as explicit full packings for $k = 2ml^2$ if $a=mb$ $(l,m \in \NN)$.
These explicit packings will give to the above quantities a
transparent geometric meaning.

The following corollary slightly refines Corollary~5.B of \cite{B1}.
\begin{corollary}   \label{c:C21}
We have \,$\max \left( 2 \dfrac{a}{b}, 8 \right) \le P(S^2(a) \times S^2(b))
\le 8 \dfrac{a}{b}$ \\
except possibly for $\frac{b}{a} = \frac{7}{8}$, in which case $P (S^2(a) \times S^2(b))
\in \{ 7,8,9 \}$.  
For $S^2(1) \times S^2(1)$  we thus have
\begin{table}[h]
 \begin{center}
 \renewcommand{\arraystretch}{1.5}
  \begin{tabular}{ccccccccc} \hline
  $k$    & 1 &      2      &     3         & 4 &      5        &     6
  & 7 & $\ge 8$      \\ \hline
  $p_k$  & $\frac{1}{2}$ & $ 1 $ & $\frac{2}{3}$ & $\frac{8}{9}$ & 
  $\frac{9}{10}$ &
  $\frac{48}{49}$ & $\frac{224}{225}$ & 1  \\ \hline
  \end{tabular}
 \end{center}
 \caption{$p_k(S^2(1) \times S^2(1))$} \label{ta:S2(1)}
\end{table}
\end{corollary}
\begin{proposition}   \label{p:C23}
Assume that $a > \frac{b}{2} > 0$.
Abbreviate 
$p_k = p_k (S^2 \ltimes S^2, \oo_{ab})$,
and set $\langle k \rangle = k$ if $k$ is odd and  $\langle k \rangle =
k+1$ if $k$ is even. Then
\[
p_k = \frac{k}{2}\frac{b}{a} \quad \mbox{ if } \; \frac{\langle k \rangle}{2}
\frac{b}{a} \le 1.
\]
Moreover,
\begin{gather*}
p_1 = \frac{b}{2a}, \qquad 
p_2 = \frac{1}{ab} \left\{b, \frac{2a+b}{4}
\right\}^2 \; \mbox{ on } \left]0,\frac{2}{3},2 \right[,      \\
p_3 = \frac{3}{2} p_2, \qquad                                            
p_4 = \frac{2}{ab} \left\{ b, \frac{2a+3b}{8}
\right\}^2 \; \mbox{ on } \left]0, \frac{2}{5}, 2 \right[,    \\ 
p_5 = \frac{5}{2ab} \left\{ b, \frac{2a+3b}{8}, \frac{2a+b}{5} 
\right\}^2 \; \mbox{ on } \left]0, \frac{2}{5}, \frac{6}{7}, 2
\right[,                                                      \\
p_6 = \frac{3}{ab} \left\{ b, \frac{2a+5b}{12},
\frac{2a+2b}{7}, \frac{2a+b}{5} 
\right\}^2 \; \mbox{ on } \left]0, \frac{2}{7}, \frac{10}{11},
\frac{4}{3}, 2  \right[,                                      \\
p_7 = \frac{7}{2ab} \left\{ b, \frac{2a+5b}{12},
\frac{6a+9b}{28}, \frac{4a+4b}{15}, \frac{4a+3b}{13}, 
\frac{6a+3b}{16}  \right\}^2                                  \\
\qquad \qquad \qquad \qquad \qquad \qquad \qquad \qquad \qquad 
\mbox{ on } 
\left]0, \frac{2}{7}, \frac{1}{2}, \frac{22}{23}, \frac{8}{7},
\frac{14}{9}, 2 \right[ .
\end{gather*}
In particular, for $k \le 7$ we have $p_k(S^2 \ltimes S^2, \oo_{ab}) =1$
exactly for 
\mbox{$(k=3, \,\frac{b}{a}=\frac{2}{3})$}, 
$(k=5, \,\frac{b}{a}=\frac{2}{5})$,
$(k=6, \,\frac{b}{a}=\frac{4}{3})$
$(k=7, \,\frac{b}{a}=\frac{2}{7})$ 
$(k=7, \,\frac{b}{a}=\frac{8}{7})$ and
\mbox{$(k=7, \,\frac{b}{a}=\frac{14}{9})$}.
\end{proposition}
In \ref{222} we will construct explicit maximal packings of $(S^2 \ltimes
S^2, \oo_{ab})$ for all $k$ with $\frac{\langle k \rangle}{2}
\frac{b}{a} \le 1$, for $k \le 5$ and $0 < \frac{b}{2} <a$ arbitrary,
and for $k=6$ and $\frac{b}{a} \in \; ]0, \frac{2}{3}] \cup [ \frac{4}{3},
2[$. Moreover, given $\oo_{ab}$ with $\frac{b}{a} = \frac{2l}{2m-l}$ for
some $l,m \in \NN$ with $m > l$, we will construct explicit full
packings of $(S^2 \ltimes S^2, \oo_{ab})$ by $l(2m-l)$ balls.
%
\begin{corollary}   \label{c:C23}
We have
$\max \left( 2 \dfrac{a}{b}, 8 \right) \le P(S^2 \ltimes S^2,
\oo_{ab})
\le  \left\{    \begin{array}{ccc}   
                       \frac{8a}{b}  & \mbox{ if } & b
                       \le a \\ [.3em]
                       \frac{8ab}{(2a-b)^2}  & \mbox{ if } & b
                       \ge a
                \end{array}   \right.         
$             
except possibly for $\frac{b}{a} \in 
\{ \frac{2}{7}, \, \frac{8}{7}, \, \frac{14}{9} \}$, in which
case the lower bound for $P \left( S^2 \ltimes S^2, \oo_{ab} \right)$
is $7$.
For $(S^2 \ltimes S^2, \oo_{11})$  we thus have
\begin{table}[h]
 \begin{center}
 \renewcommand{\arraystretch}{1.5}
  \begin{tabular}{ccccccccc} \hline
  $k$    & 1 &      2      &     3         & 4 &      5        &     6
  & 7 & $\ge 8$      \\ \hline
  $p_k$  & $\frac{1}{2}$ & $\frac{9}{16}$ & $\frac{27}{32}$ &
  $\frac{25}{32}$  &  $\frac{9}{10}$  &  $\frac{48}{49}$  &
  $\frac{14}{15}$  & 1  \\ \hline
  \end{tabular}
 \end{center}
 \caption{$p_k(S^2 \ltimes S^2, \oo_{11})$} \label{ta:ltimes12}
\end{table}
\end{corollary}
\begin{proposition}   \label{p:C24}
Let $g \ge 1$ and let $a>0$ and $b>0$. Then
\[
p_k(\Sigma_g(a) \times S^2(b)) = p_k(\Sigma_g \ltimes S^2, \oo_{ab})
= \min \left\{ 1, \frac{k}{2} \frac{b}{a} \right\}.
\]
In particular, $P(\Sigma_g (a) \times S^2(b)) = P(\Sigma_g \ltimes
S^2, \oo_{ab}) = \left\lceil \dfrac{2a}{b}  \right\rceil$.
\end{proposition}
In \ref{223} we will construct explicit maximal packings of 
$\Sigma_g(a) \times S^2(b)$ and $(\Sigma_g \ltimes S^2, \oo_{ab})$ for
all $k$ with $\left\lceil \frac{k}{2} \right\rceil \frac{b}{a} \le 1$
and explicit full packings for 
$k = 2ml^2$ if $a=mb$ or $b=ma$ $(l,m \in \NN)$.

\bigskip
\noindent
In the remainder of this section we prove Propositions~\ref{p:C21}, \ref{p:C23} and \ref{p:C24} and
Corollaries~\ref{c:C21} and \ref{c:C23}. 
We assume the reader to be familiar with \cite{B1}.
Set $N = \CC \PP^2$, let $L= [ \CC \PP^1 ]$ be the positive generator of
$H^2(N; \ZZ)$, let $\widetilde{N}_k$ be the complex blow-up of $\CC \PP^2$ at
$k$ points and let $D_1, \dots, D_k$ be the classes of the exceptional
divisors. 
\\
\\
{\it Proof of Proposition \ref{p:C21}.}
Fix $0 < b \le a$ and set $p_k = p_k(S^2(a) \times
S^2(b))$. Biran \cite[Theorem 6.1.A]{B1} showed that for any $k \in \NN$,  
\begin{equation}  \label{e:C1}
p_k = \min \left\{ 1, \frac{k}{2ab} 
\inf \left( \frac{a n_1 + b n_2}{2n_1 + 2n_2 -1} \right)^2 \right\}
\end{equation}
where the infimum is taken over all $n_1,n_2 \in \NN_0 := \NN
\cup \{ 0 \}$ for which the system
of Diophantine equations
\refstepcounter{equation}  \label{e:C1k} 
\makeatletter
\protected@edef\nummer{\theequation}
\makeatother
\begin{equation} 
\left.
  \begin{array}{rcl}   
             2 n_1 n_2       & = & \sum_{i=1}^k m_i^2 -1   \\
             2 n_1 + 2 n_2   & = & \sum_{i=1}^k m_i   +1
  \end{array}   
\right\}                     \tag{\nummer.k}
\end{equation}
has a solution $(m_1, \dots, m_k) \in \NN_0^k$.
It is easy to see that the only solutions of (\ref{e:C1k}.1) are 
($n_1=0$, $n_2=1$, $m_1=1$) and ($n_1=1$, $n_2=0$, $m_1=1$),
which yields $p_1 = \frac{b}{2a}$. 
This implies that $p_k \le \frac{k}{2} \frac{b}{a}$ for
all $0 < b \le a$ and all $k$, and that the reverse inequality holds true 
whenever $\left\lceil \frac{k}{2} \right\rceil \frac{b}{a} \le 1$
will be shown in \ref{221}.
In order to compute $p_k$ for $2 \le k \le 7$, let $\widetilde{M}_k$ be the
complex blow-up of $S^2 \times S^2$ at $k$ points and let $E_1, \dots,
E_k$ be the classes of the exceptional divisors. 
Recall that we chose the basis of $H^2(S^2 \times S^2; \ZZ)$ to be 
$\{ B = [S^2 \times pt], \, F = pt \times S^2] \}$.
The solutions of
(\ref{e:C1k}.k) correspond to the exceptional elements $E = n_1 B + n_2 F -
\sum_{i=1}^k m_i E_i \in H_2(\widetilde{M}_k; \ZZ)$ with $n_1,n_2,m_1, \dots,
m_k \ge 0$. 

Observe now that for $k \ge 1$, $\widetilde{M}_k$ is diffeomorphic to
$\widetilde{N}_{k+1}$ via a diffeomorphism under which 
the classes 
$L$, $D_1$, $D_2$, $D_3$, $\dots$, $D_{k+1}$ 
correspond to the classes 
$B+F-E_1$, $B-E_1$, $F-E_1$, $E_2$, $\dots$, $E_k$,
respectively, cf.\ Figure~\ref{figure9neu} below. 
The exceptional elements in $H_2(\widetilde{N}_{k+1}; \ZZ)$ for 
$k \le 7$ are listed in \cite[p.\ 35]{D}. 
The values $p_k$ for $2 \le k \le 7$ are now obtained by
evaluating this list in \eqref{e:C1}. 
\proofend
\\
\noindent
{\it Proof of Corollary \ref{c:C21}.}
The estimates $2 \frac{a}{b} \le P(S^2(a) \times S^2(b)) \le 8
\frac{a}{b} $ are proved in \cite[Corollary~B.5]{B1}. 
The claim now follows from the last statement in
Proposition \ref{p:C21}. 
\proofend
\\
\noindent
{\it Proof of Proposition \ref{p:C23}.}
Fix $a > \frac{b}{2} >0$ and set $M = S^2 \ltimes S^2$ and 
$p_k = p_k(M, \oo_{ab})$.
With $\aa = a - \frac{b}{2}$ and $\bb = b$ the condition $a >
\frac{b}{2} >0$ becomes $\aa >0$ and $\bb >0$.
Recall that $\oo_{ab} = bB + aF = \bb A + (\aa + \bb) F$, where 
$\{ A, F \}$ is a basis of $H_2(M; \ZZ)$ with $A \cdot A = -1$, $A \cdot F
= 1$ and $F \cdot F =0$. 
Let $\Theta \colon \widetilde{M}_k
\rightarrow M$ be the complex blow-up of $M$ at $k$ points and let $E_1,
\dots, E_k$ be the classes of the exceptional divisors. 
The first Chern class of $M$ is $c_1 = 2A + 3F$, so that the first Chern
class of $\widetilde{M}_k$ is $\tilde{c}_1 = 2A + 3F - \sum_{i=1}^k E_i$.
Let $E = n_1A+n_2F - \sum_{i=1}^k m_i E_i$ be an exceptional element in
$H_2(\widetilde{M}_k; \ZZ)$ with $m_1, \dots, m_k \ge 0$, that is,
$(n_1, n_2, m_1, \dots, m_k) \in \ZZ^2 \times \NN_0^k$ is a solution of 
the system of Diophantine equations
\refstepcounter{equation}
\label{e:C2k}
\makeatletter
\protected@edef\nummer{\theequation}
\makeatother
\begin{equation} 
\left.
  \begin{array}{rcl}   
           n_1 (2 n_2 - n_1)    & = & \sum_{i=1}^k m_i^2 -1   \\
           n_1 + 2 n_2          & = & \sum_{i=1}^k m_i   +1       
  \end{array}   
\right\} .                    \tag{\nummer.k}
\end{equation}
Suppose that $\widetilde{\oo}_{ab}$ is a symplectic form on
$\widetilde{M}_k$ such that $[\widetilde{\oo}_{ab} ] = [\Theta^*
\oo_{ab}] - \ee \sum_{i=1}^k E_i$.
We claim that for $\ee >0$ small enough, $\widetilde{\oo}_{ab}
(\tilde{c}_1 +E) >0$.
Indeed, since $m_1, \dots, m_k \ge 0$, (\ref{e:C2k}.k) implies that 
$n_1, n_2 \ge 0$. Hence
\[
\widetilde{\oo}_{ab} (\tilde{c}_1 +E)   \, =  \,         
                      \aa (2+n_1) + \bb (3+n_2) - \ee \sum_{i=1}^k (1+m_i)
\]
is positive for $\ee$ small enough.

It now follows exactly as in the proof of Theorem~6.1.A in \cite{B1}
that for any $k \in \NN$,
\begin{equation}  \label{e:C2}
p_k = \min \left\{ 1, \frac{k}{\bb (2 \aa + \bb)} 
\inf \left( \frac{\aa n_1 + \bb n_2}{n_1 + 2 n_2 -1} \right)^2 \right\}
\end{equation}
where the infimum is taken over all $n_1, n_2 \in \NN_0$ for which
(\ref{e:C2k}.k) has a solution $(m_1, \dots, m_k) \in \NN_0^k$.

Observe now that for $k \ge 0$, $\widetilde{M}_k$ is diffeomorphic to
$\widetilde{N}_{k+1}$ via a diffeomorphism under which the classes
$L$, $D_1$, $D_2$, $\dots$, $D_{k+1}$ 
correspond to the classes $A+F$, $A$, $E_1$, $\dots$, $E_k$,
respectively, see \cite[Example~7.4]{MS}.
Evaluating the list in \cite{D} in \eqref{e:C2} we obtain
\begin{gather*}
p_1 = \frac{\bb}{2\aa + \bb}, \qquad 
p_2 = \frac{2}{\bb(2\aa + \bb)} \left\{ \bb, \frac{\aa+\bb}{2}
\right\}^2 \; \mbox{ on } ]0,1,\infty[,                            \\
p_3 = \frac{3}{2} p_2, \qquad                                            
p_4 = \frac{4}{\bb(2\aa + \bb)} \left\{ \bb, \frac{\aa+2\bb}{4}
\right\}^2 \; \mbox{ on } \left]0, \frac{1}{2}, \infty \right[,    \\ 
p_5 = \frac{5}{\bb(2\aa + \bb)} \left\{ \bb, \frac{\aa+2\bb}{4}, \frac{2\aa+2\bb}{5} 
\right\}^2 \; \mbox{ on } \left]0, \frac{1}{2}, \frac{3}{2},\infty
\right[,                                                           \\
p_6 = \frac{6}{\bb(2\aa + \bb)} \left\{ \bb, \frac{\aa+3\bb}{6},
\frac{2\aa+3\bb}{7}, \frac{2\aa+2\bb}{5} 
\right\}^2 \; \mbox{ on } \left]0, \frac{1}{3}, \frac{5}{3}, 4,\infty
\right[,                                                           \\
p_7 = \frac{7}{\bb(2\aa + \bb)} \left\{ \bb, \frac{\aa+3\bb}{6},
\frac{3\aa+6\bb}{14}, \frac{4\aa+6\bb}{15}, \frac{4\aa+5\bb}{13}, 
\frac{3\aa+3\bb}{8} \right\}^2                                     \\
\qquad \qquad \qquad \qquad \qquad \qquad \qquad \qquad \qquad \qquad \mbox{ on } 
\left]0, \frac{1}{3}, \frac{2}{3}, \frac{11}{6}, \frac{8}{3}, 7, \infty \right[ .
\end{gather*}
Replacing $\aa$ by $a-\frac{b}{2}$ and $\bb$ by $b$ we finally
obtain the values $p_k$ for $1 \le k \le 7$ as stated in
Proposition~\ref{p:C23}.  
The identity $p_1 = \frac{b}{2a}$ implies that $p_k \le
\frac{k}{2} \frac{b}{a}$ for 
all $0 < \frac{b}{2} <a$ and all $k$. That the reverse inequality
holds true whenever $\frac{\langle k \rangle}{2} \frac{b}{a} \le 1$
will be shown in \ref{222}.
\proofend
\\
\noindent
{\it Proof of Corollary \ref{c:C23}.} 
Since $p_1 = \frac{b}{2a}$, we have that $P(S^2 \ltimes S^2,
\oo_{ab}) \ge \frac{2a}{b}$, and the last statement in
Proposition \ref{p:C23} shows that $P(S^2 \ltimes S^2, \oo_{ab}) \ge 8$
if $\frac{b}{a} \notin \{ \frac{2}{7}, \frac{8}{7}, \frac{14}{9} \}$ and 
$P(S^2 \ltimes S^2, \oo_{ab}) \ge 7$ if $\frac{b}{a} \in \{ \frac{2}{7}, \frac{8}{7}, \frac{14}{9} \}$.
Next, set
\[
d_{\aa \bb} = \inf \frac{\aa n_1 + \bb n_2}{n_1+2n_2-1}
\]
where the infimum is taken over all nonnegative solutions $n_1$, $n_2$,
$m_1, \dots, m_k$ of (\ref{e:C2k}.k). We claim that $d_{\aa \bb} \ge
\min \{ \aa, \frac{\bb}{2} \}$. Indeed, (\ref{e:C2k}.k) has no solution
for $n_1 = n_2 =0$. Moreover, if $m_1 = \dots = m_k =0$, then $n_1=1$ and
$n_2=0$, and the corresponding quotient is infinite. We may thus assume
that $2n_2 \ge n_1$. It is easy to see that for all $(n_1, n_2) \in \NN_0^2
\setminus \{ (0,0) \}$ with $2n_2 \ge n_1$,
\[
\frac{\aa n_1 + \bb n_2}{n_1+2n_2-1} > \min \left\{ \aa,
\frac{\bb}{2} \right\}. 
\]
Therefore, 
\[
p_k \ge \min \left\{ 1, \frac{k}{\bb (2\aa+\bb)} \min \left( \aa,
\frac{\bb}{2} \right)^2 \right\},
\]
and so
\[
P(M, \oo_{ab})
\le  \left\{    \begin{array}{ccc}   
                       \frac{4 (2 \aa + \bb)}{\bb}  & \mbox{ if } & \bb
                       \le 2 \aa \\
                       \frac{\bb (2 \aa + \bb)}{\aa^2}  & \mbox{ if } & \bb
                       \ge 2 \aa
                \end{array}   \right.      
\]    
Replacing $\aa$ by $a-\frac{b}{2}$ and $\bb$ by $b$ the claim follows.
\proofend
\\
\noindent
{\it Proof of Proposition \ref{p:C24}.} 
The statement for $\Sigma_g(a) \times S^2(b)$ was proved in
\cite[Theorem 6.1.A]{B1}.
So let $(M, \oo_{ab}) = (\Sigma_g \ltimes S^2, \oo_{ab})$. We think of 
$M$ as the projectivization 
$\PP (L_1 \oplus \CC) \xrightarrow{\pi} \Sigma_g$ 
of the complex rank two bundle $L_1 \oplus \CC$ over $\Sigma_g$, where
$L_1$ is a holomorphic line bundle of Chern index $1$,
and we endow $\PP (L_1 \oplus \CC)$ with its
canonical complex structure $J_{can}$.
Let $(\widetilde{M}_k, \widetilde{J}_{can})$ be the complex blow-up of
$(\PP (L_1 \oplus \CC), J_{can})$ at $k$ generic points and let
$\widetilde{\oo}_{ab}$ be a blow-up of $\oo_{ab}$. 
Finally, denote by $\ce_k$ the set of homology classes of $\widetilde{M}_k$
which can be represented by $\widetilde{\oo}_{ab}$-symplectic 
exceptional spheres. 
We claim that
\[
\ce_k (\widetilde{M}_k, \widetilde{\oo}_{ab}) = \{ E_1, \dots, E_k,
F-E_1, \dots, F-E_k \},
\]
where $E_1, \dots, E_k$ are the classes of the exceptional divisors and 
$F \in H^2(M; \ZZ) \subset H^2(\widetilde{M}_k; \ZZ)$ is the fibre
class. 
The analogous statement for $\Sigma_g(a) \times S^2(b)$ was proved by
Biran, \cite{B1}, 
in the proof of his Corollary~5.C. His argument immediately
applies to the twisted bundle and is repeated here for the sake of its beauty.

So let $E= n_1 A + n_2 F - \sum_{i=1}^k m_i E_i \in \ce_k$. 
Let $\cj (\widetilde{\oo}_{ab})$ be the space of 
$\widetilde{\oo}_{ab}$-tamed almost complex structures on 
$\widetilde{M}_k$ and let $\cj_E
\subset \cj (\widetilde{\oo}_{ab})$ be the subset of those  $J$ for
which there exist $J$-holomorphic $E$-spheres. It is known \cite[Chapter
V, proof of Lemma 2.C.2]{B0} that $\cj_E$ contains a path-connected set
which is open and dense in $\cj (\widetilde{\oo}_{ab})$. 
Let $\{ \widetilde{J}_t \}_{0 \le t \le 1}$ be a smooth path in $\cj
(\widetilde{\oo}_{ab})$ with $\{ \widetilde{J}_t \}_{0 \le t < 1}
\subset \cj_E$ and $\widetilde{J}_1 = \widetilde{J}_{can}$. 
Gromov's compactness theorem now shows that there exists a connected
(but possibly cusp) $\widetilde{J}_{can}$-holomorphic $E$-curve $C = C_1
\cup \dots \cup C_n$ with $g(C_j) =0$ for all $j$. 

Let $\tilde{\pi} \colon \widetilde{M}_k \ra \Sigma_g$ be the lift of
$\pi \colon M \ra \Sigma_g$. Since $\pi$ is $J_{can}$-holomorphic,
$\tilde{\pi}$ is $\widetilde{J}_{can}$-holomorphic.
Let $h_j \colon S^2 \ra C_j$ be a $\widetilde{J}_{can}$-holomorphic
parametrization of $C_j$ and let $l_j \colon S^2 \ra \CC$ be a lift of
$\tilde{\pi} \circ h_j$ to the universal cover of $\Sigma_g$. By
Liouville's theorem, $l_j$ is constant, and so $\tilde{\pi}(C_j)$ is a
point. Since this holds true for all $j$ and since $C$ is connected,
$\tilde{\pi}(C)$ is a point too. Hence, and since $A$ is the class of a
section, 
\[
0 = \tilde{\pi}_* \left( [C] \right) = \tilde{\pi}_* (E) = \pi_* (n_1 
A + n_2 F) = n_1 [\Sigma_g],
\]
and so $n_1 =0$, i.e., $E= n_2 F - \sum_{i=1}^k m_i E_i$. 
Since the first Chern class of $\widetilde{M}_k$ is $\tilde{c}_1 = 2A +
(3-2g)F-\sum_{i=1}^k E_i$, the conditions $E \cdot E = -1$ and $\tilde{c}_1(E)
=1$ become $\sum_{i=1}^k m_i^2=1$ and $2n_2 - \sum_{i=1}^k m_i =1$,
which implies $E \in \{ E_1, \dots, E_k, F-E_1, \dots, F-E_k \}$.

Conversely, $E_i$ is clearly an $\widetilde{\oo}_{ab}$-symplectic
exceptional class, and the proper
transform of the $J_{can}$-fibre passing through the point $\Theta_* (E_i)$ is a
$\widetilde{J}_{can}$-exceptional rational curve and hence an
$\widetilde{\oo}_{ab}$-symplectic exceptional sphere in class $F-E_i$.

Finally, we have that $\oo_{ab}(F) = b$, $c_1(F) =2$ and 
$2 \Vol (M, \oo_{ab}) = 2ab$.   
Proposition~\ref{p:C24} now follows from Theorem~6.A in \cite{B1}.
\proofend

\section{Explicit maximal packings in four dimensions}  \label{s9:3}
\label{s:construction}

\ni
In this section we realize most of the packing numbers computed
in the previous section by explicit symplectic packings. 
Sometimes, we shall give two different maximal packings. It is
known that for the $4$-ball, $\CC \PP^2$ and ruled symplectic $4$-manifolds,
any two packings by $k$ balls of equal size are symplectically isotopic,
see \cite{B-96,M2}.

Recall that $\RR^4$ is endowed with the symplectic form $\oo_0 =
dx_1 \wedge dy_1 + dx_2 \wedge dy_2$.
We shall often use the Lagrangian splitting $\RR^2(x) \times
\RR^2(y)$ of $\RR^4$.
Set $\sq^2(1) = \;]0,1[\, \times\, ]0,1[\; \subset \RR^2(y)$. 
In order to construct our symplectic packings by balls, we shall
construct explicit symplectic embeddings of a ball $B^4(a)$ into
products $U \times \sq^2(1)$ of almost equal volume, where $U
\subset \RR^2(x)$ is a domain as in Figure~\ref{figure41} below.
The symplectic $4$-manifolds $(M,\oo)$ we shall consider contain a
domain of equal volume which is explicitly symplectomorphic to $V
\times \sq^2(1) \subset \RR^2(x) \times \RR^2(y)$.
In order to construct an explicit symplectic packing of $(M,\oo)$ by $k$
equal balls it will thus suffice to insert $k$ disjoint domains
$U$ of equal width as in Figure~\ref{figure41} into $V$.

In the explicit packings constructed in \cite{Ka,T,Kr,MMT}, a
ball is viewed as a product $U \times \sq^2(1)$, where $U$ is an
affine image of a simplex and thus in particular convex.
Our domains $U$ need not be convex, and so we have a larger
arsenal of shapes at our disposal.

\subsection{How to map $B^4(a)$ to $U \times \sq^2(1)$.}
\label{s:construction.1}

Let $D(a)$ be the open disc in $\RR^2$ of area $a$ centred at
the origin, and let 
\[
R(a) \,=\, \left\{ (x,y) \in \RR^2 \mid 0 < x<a,\, 0<y<1 \right\}.
\]
Our symplectic embeddings $B^4(a) \ha U \times \sq^2(1) \subset
\RR^2(x) \times \RR^2(y)$ will be obtained from restricting split
symplectomorphisms $\aa_1 \times \aa_2 \colon D(a) \times D(a)
\ra R(a) \times R(a)$ to $B^4(a)$.
Notice that in dimension $2$ an embedding is symplectic if and
only if it is area and orientation preserving.
In order to explicitly describe such embeddings, we follow \cite{S-Israel}
and start with

\begin{definition}
{\rm
A family $\cl$ of loops in $R(a)$ is {\it admissible}\, if there exists
a diffeomorphism $\bb \colon D(a) \setminus \{0\} \ra R(a) \setminus \{ p
\}$ for some point $p \in R(a)$ such that 
\begin{itemize}
\item[(i)]
concentric circles are mapped to elements of $\cl$,
\item[(ii)]
in a neighbourhood of the origin $\bb$ is a translation.
\end{itemize}
}
\end{definition}
\begin{lemma} \label{lemmaarea}
\index{area preserving embedding!construction of|(}
Given an admissible family $\cl$ of loops in $R(a)$,
there exists a symplectomorphism from $D(a)$ to $R(a)$
mapping concentric circles to loops of $\cl$.
\end{lemma}
We refer to \cite{S-Israel} or \cite[Lemma~2.5]{buch} for the
elementary proof.
Notice that the symplectomorphism guaranteed by the lemma is
uniquely determined by its image of the ray 
$\left\{ (x,0) \in D(a) \mid x \ge 0 \right\}$. 
We can thus ``explicitly'' describe a symplectomorphism from
$D(a)$ to $R(a)$ by prescribing an admissible family of loops in
$R(a)$ and a smooth line from the centre of $\cl$ to the
boundary of $R(a)$ meeting each loop exactly once.              
\begin{figure}[h] 
       \begin{center}
  \psfrag{a}{$(\aa)$}
  \psfrag{aa}{$(\overline{\aa})$}
  \psfrag{b}{$(\bb)$}
  \psfrag{c}{$(\gg)$}
  \psfrag{d}{$(\dd)$}
  \leavevmode\epsfbox{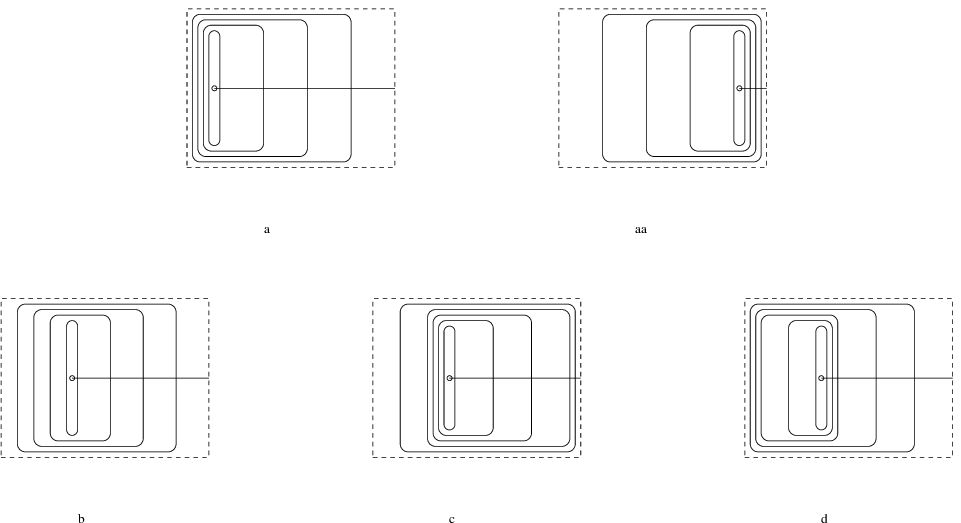}
 \end{center}
 \caption{Symplectomorphisms $D(a) \ra R(a)$.} \label{figure40}
\end{figure}
%
%
For the symplectomorphisms $D(a) \ra R(a)$ described in Figure~\ref{figure40}
we have chosen this line to be a segment at height $y= \frac 12$
from the centre of $\cl$ to the right boundary.
Consider first the symplectomorphism $\aa$ represented by
$(\aa)$ in Figure~\ref{figure40}.
The restriction of $\aa \times \aa$ to $B^4(a)$ is contained in
the product $U \times \sq^2(1) \subset \RR^2(x) \times
\RR^2(y)$, where $U$ is a small neighbourhood of the simplex $[
\aa \aa ]$ of width $a$ shown in Figure~\ref{figure41}.
Since for every neighbourhood $U$ of $[ \aa \aa ]$ we can choose
$\aa$ such that $( \aa \times \aa ) \left( B^4(a) \right)$ is
contained in $U \times \sq^2(1)$, we shall work with the
simplex $[\aa \aa ]$ instead of $U$.
The bar in the notation $\overline{\aa}$ used in Figure~\ref{figure40}
and Figure~\ref{figure41} indicates that $\overline{\aa}$ is
the mirror of $\aa$.
Figure~\ref{figure41} shows the $x_1$-$x_2$-shadows of the image of
$B^4(a)$ under some other products of the symplectomorphisms in
Figure~\ref{figure40} and of their mirrors. We invite the reader
to create further shadows.

\begin{figure}[h] 
 \begin{center}
  \psfrag{a}{$[ \aa \aa ]$}
  \psfrag{aa}{$[\overline{\aa} \aa ]$}
  \psfrag{aaa}{$[ \overline{\aa} \,\overline{\aa} ]$}
  \psfrag{b}{$[ \bb \bb ]$}
  \psfrag{ca}{$[ \gg \aa ]$}
  \psfrag{dd}{$[ \overline{\dd} \dd ]$}
  \leavevmode\epsfbox{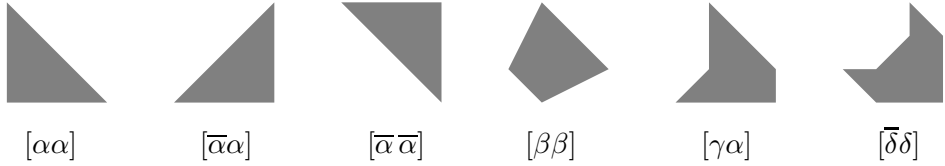}
 \end{center}
 \caption{Some $x_1$-$x_2$-shadows.} \label{figure41}
\end{figure}
%
%

\begin{remark}  \label{r9:explicit}
{\rm
Besides of being explicit, the $4$-dimensional symplectic packings
constructed in \cite{Ka,T} and in this section have yet another
advantage over the packings found in \cite{MP,B1,B2}:
The symplectic packings of $(M,\oo)$ by $k$ balls obtained from the
method in \cite{MP,B1,B2} are maximal in the following sense.
For every $\ee >0$ there exists a symplectic embedding 
$\ff_\ee \colon \coprod_{i=1}^k B^{2n}(a) \ha (M,\oo)$ such that 
\begin{equation}  \label{e9:full}
\frac{\Vol \left( \im \ff_\ee, \oo \right)}{\Vol \left(M,\oo \right)}
\,\ge\, p_k(M,\oo) - \ee .
\end{equation}
Karshon's symplectic packings of $\left( \CC \PP^2,\oo_{SF} \right)$ 
by $2$ and $3$ balls $B^4 \left( \frac{\pi}{2}
\right)$ given by the map \eqref{emb:BCP2} and compositions of this map with
coordinate permutations fill {\it exactly}\, $\frac 12$ and $\frac
34$ of $\left( \CC \PP^2, \oo_{SF} \right)$.
Similarly, the 4-dimensional packings in \cite{T} and in this section are
maximal in the following sense:

\m
\ni
{\it
There exists a symplectic embedding $\ff \colon \coprod_{i=1}^k B^4(a)
\ha (M,\oo)$ such that 
\begin{equation}  \label{e9:fullee}
\frac{\Vol \left( \im \ff, \oo \right)}{\Vol \left(M,\oo \right)}
\,=\, p_k(M,\oo) .
\end{equation}
Moreover, $\ff$ is explicit in the following sense: 
The image $\coprod_{i=1}^k \ff \left( B^4(a) \right)$ of $\ff$ is
explicit, and given $a' < a$ one can construct $\ff$ such that its
restriction to $\coprod_{i=1}^k B^4(a')$ is given pointwise.
}

\m
\ni
Indeed, choose a sequence $a'<a_j \nearrow a$. The packings in \cite{T}
and our packings $\ff(a_j) \colon \coprod_{i=1}^k B^4(a_j) \ha (M,\oo)$
can be chosen such that
\[
\im \ff (a_j) \subset \im \ff (a_{j+1})
\quad\, \text{ for all } \, j.
\]
The claim now follows from a result of Mc\:\!Duff, \cite{M}, 
stating that two symplectic embeddings of a closed ball into a larger ball are isotopic
via a symplectic isotopy of the larger ball.
\diam
}
\end{remark}


%
%
%
%


\subsection{Maximal packings of the $4$-ball and of $\CC \PP^2$}  \label{21}

\begin{figure}[h] 
 \begin{center}
  \psfrag{a}{$(a)$}
  \psfrag{b}{$(b)$}
  \psfrag{c}{$(c)$}
  \leavevmode\epsfbox{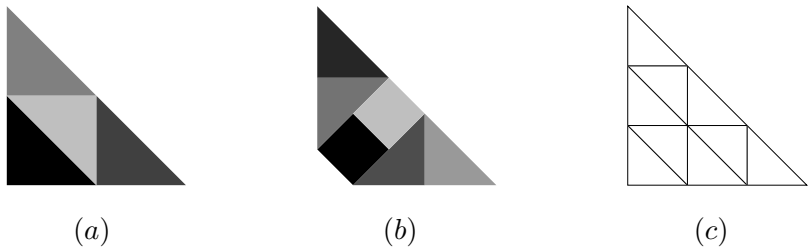}
 \end{center}
 \caption{Maximal packings of $B^4$ for $k \le 6$ and $k = l^2$.} \label{figure42}
\end{figure}

In view of the symplectomorphism \eqref{emb:BCP2} and the
identity \eqref{pk=} we only need to construct packings of the
$4$-ball.
It follows from Table~\ref{ta:ball} that any $k$ of the embeddings in
\mbox{Figure \ref{figure42}(a)} yield a maximal packing of 
$B^4$ by $k$ balls,
$k = 2,3,4$, and that any $k$ of the embeddings in \mbox{Figure
\ref{figure42}(b)} yield a maximal packing by $k = 5,6$ balls. 
\mbox{Figure~\ref{figure42}(c)} shows a full packing by 9 balls.
%

Explicit maximal packings of $B^4$ by $k \le 6$ balls were first
constructed by Traynor in \cite{T}.
Her packings by $5$ or $6$ balls are constructed by a Lagrangian
folding method.
Neither Traynor's nor our packing method nor their combination
can realize the packing numbers 
$p_7 \left( B^4 \right) = \frac{63}{64}$ and $p_8 \left( B^4
\right) = \frac{288}{289}$, but they   
only fill $\frac 79$ and $\frac 89$ of the
$4$-ball by $7$ and $8$ equal balls, respectively. 

\begin{question}  \label{q:78}
Is there an explicit embedding of $7$ or $8$ equal balls into
the $4$-ball filling more than $\frac 79$ and $\frac 89$ of the volume?
\end{question}

\subsection{Maximal packings of ruled symplectic $4$-manifolds}
\label{22}

Given a ruled symplectic $4$-manifold $(M, \oo_{ab})$, let $c_k(a,b)$ be the
supremum of those $A$ for which $\coprod_{i=1}^k B^{2n}(A)$
symplectically embeds into $(M,\oo_{ab})$, so that
\begin{equation}  \label{e:pkc}
p_k(M, \oo_{ab}) = \frac{k \, c_k^2(a,b)}{2 \Vol (M, \oo_{ab})} .
\end{equation}
We shall write $c$ instead of $c_k(a,b)$ if
$\left(M,\oo_{ab}\right)$ and $k$ are clear from the context.

\subsubsection{Maximal packings of $S^2(a) \times S^2(b)$.}
\label{221}

As in Proposition~\ref{p:C21} we assume that $a \ge b$.
Represent the symplectic structure of $S^2(a) \times S^2(b)$ by a split form.
Using Lemma~\ref{lemmaarea} we symplectically identify 
$S^2(a) \setminus pt$ with $]0,a [ \times ]0,1[$ and 
$S^2(b) \setminus pt$ with $]0,b [ \times ]0,1[$. Then
\[
\sq (a,b) \times \sq^2(1) \,=\, 
S^2(a) \times S^2(b) \setminus \left\{ S^2(a) \times pt \cup pt
\times S^2(b) \right\} .
\] 
Besides for $k \in \{ 6,7 \}$, we will construct the explicit
maximal packings promised after Proposition~\ref{p:C21} by 
constructing packings of $\sq (a,b) \times \sq^2(1)$
which realize the packing numbers of $S^2(a) \times S^2(b)$ computed in 
Proposition~\ref{p:C21} and hence are maximal. 
(It is, in fact, known that {\it all}\, packing numbers of 
$\sq (a,b) \times \sq^2(1)$ and $S^2(a) \times S^2(b)$ agree, 
see \cite[Remark 2.1.E]{MP}).

\begin{figure}[h] 
 \begin{center}
  \psfrag{a}{$(a)$}
  \psfrag{b}{$(b)$}
  \psfrag{aa}{$a$}
  \psfrag{bb}{$b$}
  \psfrag{x1}{$x_1$}
  \psfrag{x2}{$x_2$}
  \leavevmode\epsfbox{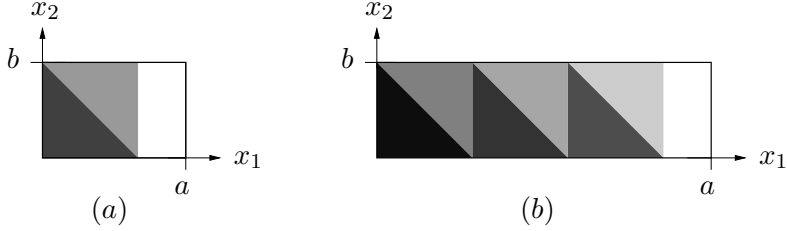}
 \end{center}
 \caption{Maximal packings of $S^2(a) \times S^2(b)$ by $k$ balls, 
          $\left\lceil \frac{k}{2} \right\rceil \frac{b}{a} \le 1$.} \label{figure43}
\end{figure}

\s
To construct explicit maximal packings for all $k$ with $\left\lceil \frac{k}{2} \right\rceil
\frac{b}{a} \le 1$ is a trivial matter. \mbox{Figure \ref{figure43}} shows a
maximal packing by $1$ and $2$ respectively 5 and 6 balls.
%
%

\begin{figure}[h] 
 \begin{center}
  \psfrag{x1}{$x_1$}
  \psfrag{x2}{$x_2$}
  \psfrag{aa}{$a$}
  \psfrag{bb}{$b$}
  \psfrag{12}{$\frac{b}{a} = \frac{1}{2}$}
  \psfrag{34}{$\frac{b}{a} = \frac{3}{4}$}
  \psfrag{11}{$\frac{b}{a} = 1$}
  \leavevmode\epsfbox{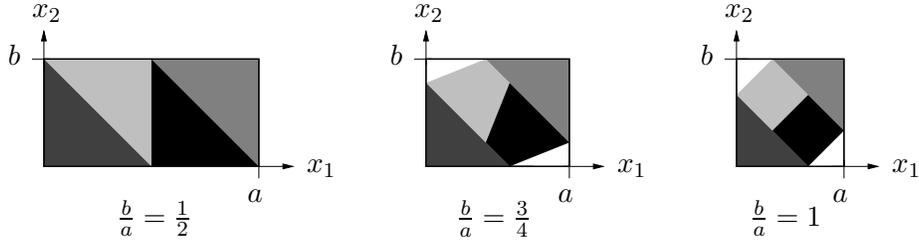}
 \end{center}
 \caption{Maximal packings of $S^2(a) \times S^2(b)$ by 3 and 4 balls, 
          $\frac{b}{a} \ge \frac{1}{2}$.} \label{figure44}
\end{figure}

Let now $k = 3,4$ and $\frac{b}{a} \ge \frac{1}{2}$. \mbox{Figure
\ref{figure44}} shows maximal packings of $S^2(a) \times S^2(b)$ by $k$ balls
for $\frac{b}{a} = \frac{1}{2}$, $\frac{b}{a} = \frac{3}{4}$ and
$\frac{b}{a} = 1$.
%
%
For $\frac{b}{a} > \frac{1}{2}$ the $(x_1, x_2)$-coordinates of the
vertices of the ``upper left ball'' are 
\[
(0, c), \quad (a-c, b),\quad (c, c),\quad (a-c, b-c),
\]
where $c = \frac{a+b}{3}$.
As in most of the subsequent figures, 
the three pictures in Figure~\ref{figure44} should be seen as
moments of a movie  
starting at $\frac ba = \frac 12$ and ending at $\frac ba =
1$. Each ball in this movie moves in a smooth way.

\begin{figure}[h] 
 \begin{center}
  \psfrag{x1}{$x_1$}
  \psfrag{x2}{$x_2$}
  \psfrag{aa}{$a$}
  \psfrag{bb}{$b$}
  \psfrag{13}{$\frac{b}{a} = \frac{1}{3}$}
  \psfrag{12}{$\frac{b}{a} = \frac{1}{2}$}
  \psfrag{34}{$\frac{b}{a} = \frac{3}{4}$}
  \psfrag{11}{$\frac{b}{a} = 1$}
  \leavevmode\epsfbox{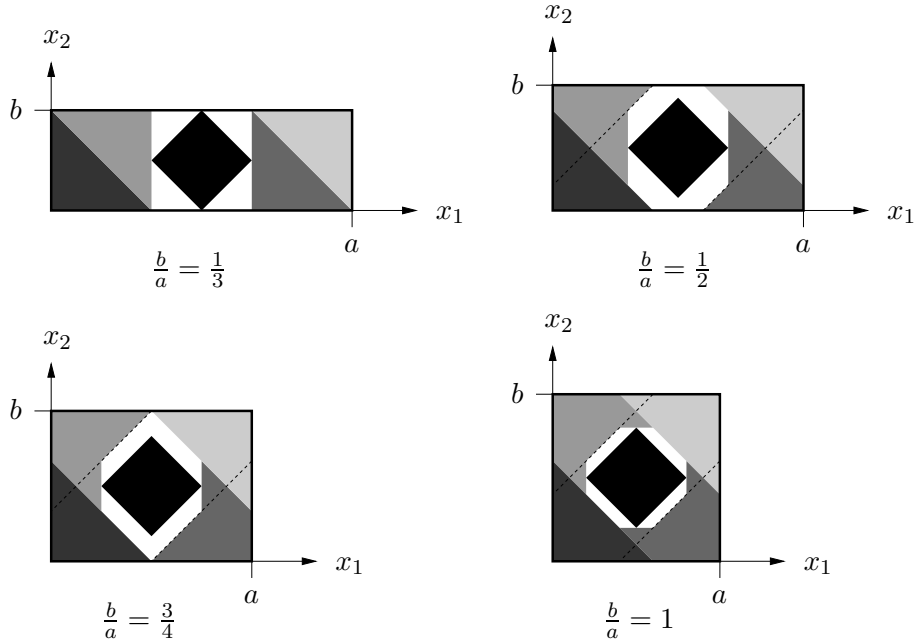}
 \end{center}
 \caption{Maximal packings of $S^2(a) \times S^2(b)$ by 5 balls, 
          $\frac{b}{a} \ge \frac{1}{3}$.} \label{figure45}
\end{figure}

Next, let $k=5$ and $\frac{b}{a} \ge \frac{1}{3}$. 
In order to construct a smooth family of maximal packings of $S^2(a)
\times S^2(b)$ by 5 balls, we think of the maximal packing for
$\frac{b}{a} = \frac{1}{3}$ rather as in Figure~\ref{figure45}
than as in Figure~\ref{figure43}(a). 
The $x_1$-width of all balls is $\frac{a+ 2b}{5}$, 
and the ``upper left ball'' has 5 vertices for 
$\frac{1}{3} < \frac{b}{a} \le \frac{3}{4}$ and $7$ vertices for 
$\frac{b}{a} > \frac{3}{4}$.
%

For $k \in \{ 6,7 \}$, we cannot realize the packing numbers
$p_k \left( S^2(a) \times S^2(b) \right)$ by directly packing
rectangles as for $k \le 4$.
We shall instead construct certain maximal packings of $\CC
\PP^2$ which correspond to maximal packings of $S^2(a) \times
S^2(b)$.
As noticed in \cite{B1}, the correspondence between symplectic
packings and the symplectic blow-up operation
%
and the diffeomorphism mentioned in the proof of Proposition~\ref{p:C21}
imply
\begin{lemma}  \label{lemma:SSCP}
Packing $S^2(a) \times S^2(b)$ by $k$ equal balls
$\coprod_{i=1}^k B^4(c)$ corresponds to packing 
$\left( \CC \PP^2, (a+b-c) \oo_{SF} \right)$ by the $k+1$ balls 
$B^4(a-c) \coprod B^4(b-c) \coprod_{i=1}^{k-1} B^4(c)$.
\end{lemma}

\begin{figure}[h] 
 \begin{center}
  \psfrag{x1}{$x_1$}
  \psfrag{x2}{$x_2$}
  \psfrag{a}{$a$}
  \psfrag{b}{$b$}
  \psfrag{c}{$c$}
  \psfrag{abc}{$a+b-c$}
  \psfrag{ac}{$a-c$}
  \psfrag{bc}{$b-c$}
    \leavevmode\epsfbox{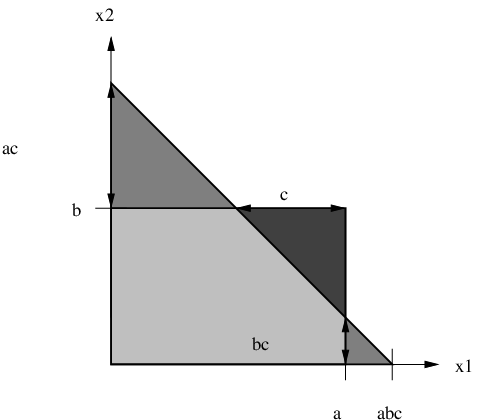}
 \end{center}
 \caption{$\left( \CC \PP^2, (a+b-c) \oo_{SF} \right) \setminus
           B^4(a-c) \coprod B^4(b-c) \,=\, 
S^2(a) \times S^2(b) \setminus B^4(c)$} \label{figure9neu}
\end{figure}

In order to make this correspondence plausible, we choose $\frac
ba = \frac 23$ and $c=c_6(a,b) = \frac{a+2b}{5}$, and we think
of $\left( \CC \PP^2, (a+b-c) \oo_{SF} \right)$ as the simplex
of width $a+b-c$ and of $S^2(a) \times S^2(b)$ as the rectangle
of width $a$ and length $b$. As Figure~\ref{figure9neu}
illustrates, the space obtained by removing a ball $B^4(c)$ from
$S^2(a) \times S^2(b)$ coincides with the space obtained by
removing the balls
$B^4(a-c) \coprod B^4(b-c)$ from $\left( \CC \PP^2, (a+b-c)
\oo_{SF} \right)$.
%
%

\begin{figure}[h] 
 \begin{center}
  \psfrag{x1}{$x_1$}
  \psfrag{x2}{$x_2$}
  \psfrag{abc}{$a+b-c$}
  \psfrag{ac}{$a-c$}
  \psfrag{13}{$\frac{b}{a} = \frac{1}{3}$}
  \psfrag{12}{$\frac{b}{a} = \frac{1}{2}$}
  \psfrag{23}{$\frac{b}{a} = \frac{2}{3}$}
  \leavevmode\epsfbox{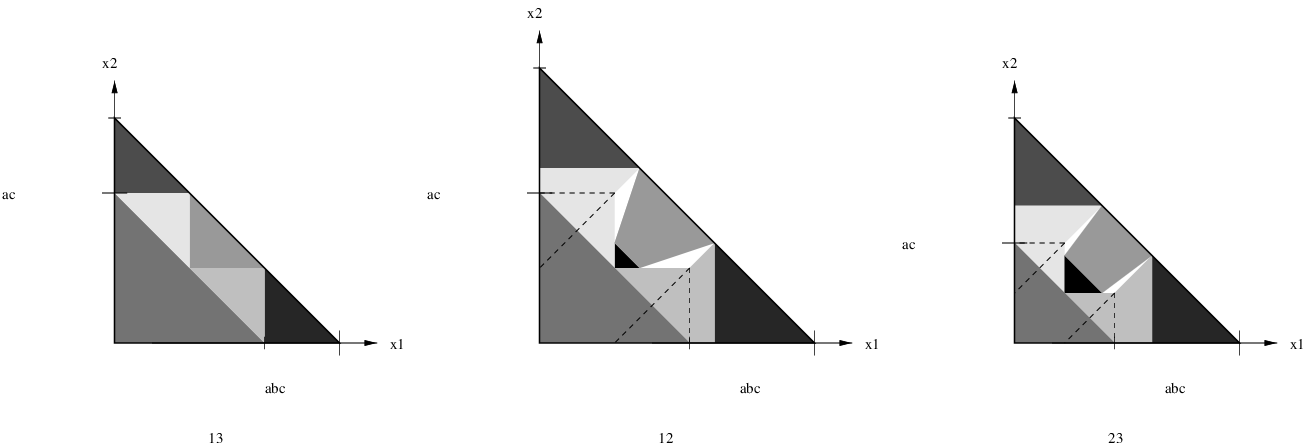}
 \end{center}
 \caption{Maximal packings of $S^2(a) \times S^2(b)$ by $6$ balls, 
          $\frac{1}{3} \le \frac{b}{a} \le \frac{3}{4}$.} \label{figure80}
\end{figure}

\begin{figure}[h] 
 \begin{center}
  \psfrag{x1}{$x_1$}
  \psfrag{x2}{$x_2$}
  \psfrag{abc}{$a+b-c$}
  \psfrag{ac}{$a-c$}
  \psfrag{34}{$\frac{b}{a} = \frac{3}{4}$}
  \psfrag{45}{$\frac{b}{a} = \frac{4}{5}$}
  \psfrag{11}{$\frac{b}{a} = 1$}
  \leavevmode\epsfbox{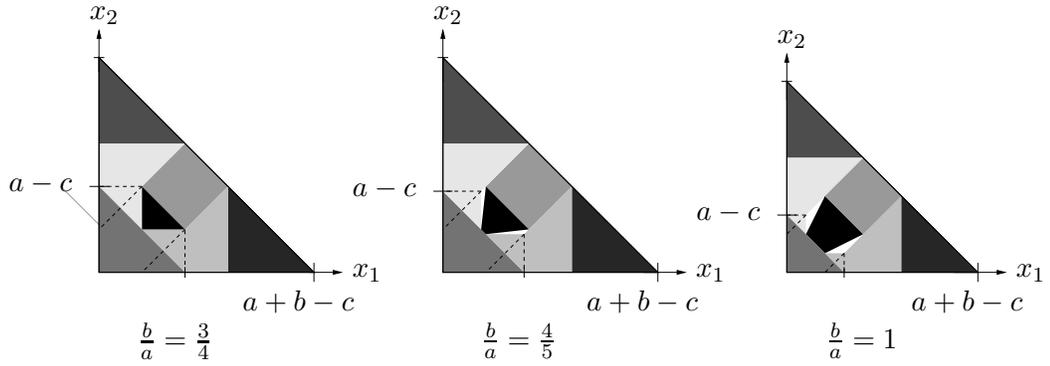}
 \end{center}
 \caption{Maximal packings of $S^2(a) \times S^2(b)$ by $6$ balls, 
          $\frac{3}{4} \le \frac{b}{a} \le 1$.} \label{figure81}
\end{figure}

\begin{figure}[h] 
 \begin{center}
  \psfrag{x1}{$x_1$}
  \psfrag{x2}{$x_2$}
  \psfrag{abc}{$a+b-c$}
  \psfrag{ac}{$a-c$}
  \psfrag{14}{$\frac{b}{a} = \frac{1}{4}$}
  \psfrag{25}{$\frac{b}{a} = \frac{2}{5}$}
  \psfrag{35}{$\frac{b}{a} = \frac{3}{5}$}
  \leavevmode\epsfbox{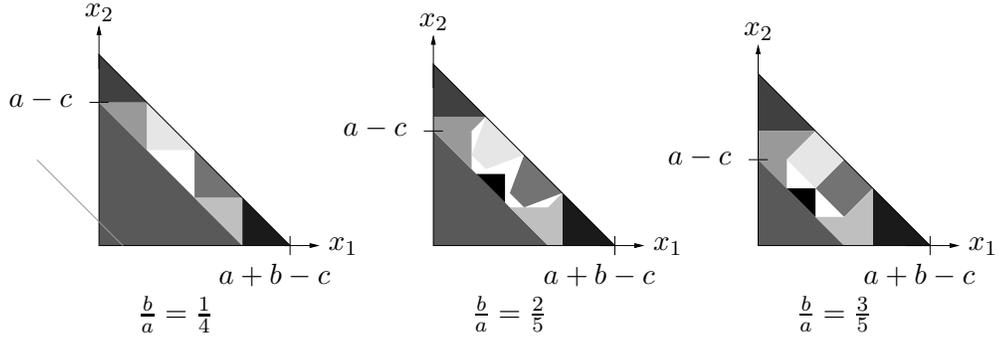}
 \end{center}
 \caption{Maximal packings of $S^2(a) \times S^2(b)$ by $7$ balls, 
          $\frac{1}{4} \le \frac{b}{a} \le \frac{3}{5}$.} \label{figure777}
\end{figure}

Figures~\ref{figure80}, \ref{figure81} and \ref{figure777}
describe explicit packings of 
$\left( \CC \PP^2, (a+b-c) \oo_{SF} \right)$ by balls 
$B^4(a-c) \coprod B^4(b-c) \coprod_{i=1}^{k-1} B^4(c)$ for $k
\in \{6,7 \}$ and $c$ as in Proposition~\ref{p:C21}.
The lower left triangle represents $B^4(a-c)$ and the black "ball"
represents $B^4(b-c)$. 
%
%
%
%
%

From these packings one obtains explicit packings of $S^2(a)
\times S^2(b)$ as follows:
First symplectically blow up $\left( \CC \PP^2, (a+b-c) \oo_{SF}
\right)$ twice by removing the balls $B^4(a-c)$ and $B^4(b-c)$
and collapsing the remaining boundary spheres to exceptional
spheres in homology classes $D_1$ and $D_2$. 
The resulting manifold, which is symplectomorphic to 
$S^2(a) \times S^2(b)$ blown up at one point with weight $c$,
still contains the $k-1$ explicitly embedded balls $B^4(c)$,
and according to \cite[Theorem~4.1.A]{B1} the exceptional
sphere in class $L-D_1-D_2$ can be symplectically blown down with weight
$c$ to yield the $k$'th ball $B^4(c)$ in $S^2(a) \times S^2(b)$.
%

%
%
%

\begin{figure}[h] 
 \begin{center}
  \psfrag{x1}{$x_1$}
  \psfrag{x2}{$x_2$}
  \psfrag{b}{$b$}
  \psfrag{2b}{$2b$}
  \leavevmode\epsfbox{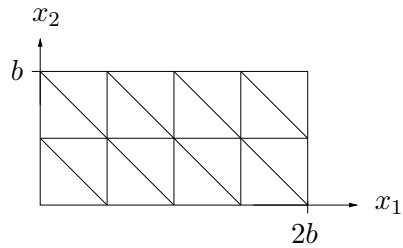}
 \end{center}
 \caption{A full packing of $S^2(2b) \times S^2(b)$ by 16 balls.}
 \label{figure46} 
\end{figure}

Finally, the construction of full packings of $S^2(mb) \times S^2(b)$ by
$2ml^2$ balls $(l,m \in \NN)$ is also straightforward. 
Figure~\ref{figure46} shows such a packing for $l=m=2$. 
%
%

\subsubsection{Maximal packings of $\left(S^2 \ltimes S^2, \oo_{ab}\right)$.}
\label{222}

In order to describe our maximal packings of 
$\left(S^2 \ltimes S^2, \oo_{ab}\right)$, it will be convenient
to work with the parameters 
$\aa = a - \frac b2$, $\bb = b$, so that $\aa>0$, $\bb>0$ and
$\oo_{ab} = \bb A + (\aa+\bb)F$. 
Recall that $S^2 \ltimes S^2$ is diffeomorphic to the 
blow-up $\widetilde{N}_1$ of $\CC \PP^2$ at one point 
via a diffeomorphism under which 
$L$, $D_1$ correspond to $A+F$, $A$. 
We can therefore view $\left(S^2 \ltimes S^2, \oo_{ab}\right)$
as $\widetilde{N}_1$ endowed with the symplectic form in class 
$(\aa+\bb)L - \aa D_1$ obtained by symplectically blowing up $\left( \CC
\PP^2, (\aa + \bb) \oo_{SF} \right)$ with weight $\aa$.
Since symplectic blowing up with weight $\aa$ corresponds to
removing a ball $B^4(\aa)$ and collapsing the remaining boundary
sphere to an exceptional sphere in class $D_1$, 
we can think of this symplectic manifold as the truncated
simplex obtained by removing the simplex of width $\aa$ from the
simplex of width $\aa +\bb$.
  
Denote by $\lfloor x \rfloor$ the integer part of $x \ge 0$.
In the parameters $\aa$ and $\bb$, the packings promised  
after Proposition~\ref{p:C23} are explicit
maximal packings of $(S^2 \ltimes S^2, \oo_{ab})$ for all $k$ with
$\lfloor \frac{k}{2} \rfloor \frac{\bb}{\aa} \le 1$, 
for $k \le 5$ and $\aa, \bb
>0$ arbitrary, and for $k=6$ and $\frac{\bb}{\aa} \in \; ]0,1] \cup [4, \infty
[$. Moreover, given $\oo_{ab}$ with $\frac{\bb}{\aa} = \frac{l}{m-l}$ for
some $l,m \in \NN$ with $m > l$, we will construct explicit full
packings of $(S^2 \ltimes S^2, \oo_{ab})$ by $l(2m-l)$ balls.

Set $c_k = c_k(a,b) = c_k(S^2 \ltimes S^2, \oo_{ab})$.
Using $2 \Vol \left( S^2 \ltimes S^2, \oo_{ab} \right) = \bb
(2\aa+\bb)$ and \eqref{e:pkc} we read off from the list in 
the proof of Proposition \ref{p:C23} that
\begin{gather*}
c_1 = \bb, \qquad 
c_2 = c_3 = \left\{ \bb, \frac{\aa+\bb}{2} \right\} \; \mbox{ on }
]0,1,\infty[,   \\    
c_4 = \left\{ \bb, \frac{\aa+2\bb}{4} \right\} \; \mbox{ on } \left]0, \frac{1}{2}, \infty \right[, \\ 
c_5 = \left\{ \bb, \frac{\aa+2\bb}{4}, \frac{2\aa+2\bb}{5} 
\right\} \; \mbox{ on } \left]0, \frac{1}{2}, \frac{3}{2},\infty \right[,                           \\
c_6 = \left\{ \bb, \frac{\aa+3\bb}{6}, \frac{2\aa+3\bb}{7}, \frac{2\aa+2\bb}{5} 
\right\} \; \mbox{ on } \left]0, \frac{1}{3}, \frac{5}{3}, 4,\infty
\right[.
%
\end{gather*}

\begin{figure}[h] 
 \begin{center}
  \psfrag{x1}{$x_1$}
  \psfrag{x2}{$x_2$}
  \psfrag{aa}{$\aa$}
  \psfrag{ab}{$\aa + \bb$}
  \psfrag{a}{$(a)$}
  \psfrag{b}{$(b\,1)$}
  \psfrag{c}{$(b\,2)$}
  \leavevmode\epsfbox{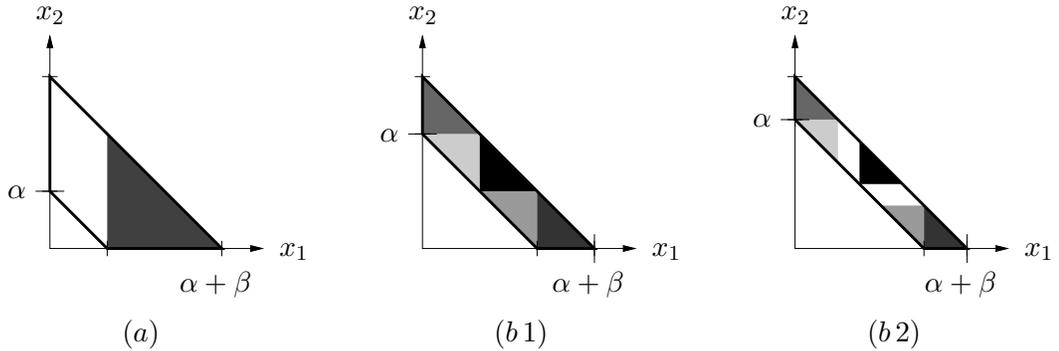}
 \end{center}
 \caption{Maximal packings of $(S^2 \ltimes S^2, \oo_{ab})$ by $k$
          balls, $\lfloor \frac{k}{2} \rfloor \frac{\bb}{\aa} \le 1$.}
                                   \label{figure47}
\end{figure}

\begin{figure}[h] 
 \begin{center}
  \psfrag{x1}{$x_1$}
  \psfrag{x2}{$x_2$}
  \psfrag{aa}{$\aa$}
  \psfrag{ab}{$\aa + \bb$}
  \leavevmode\epsfbox{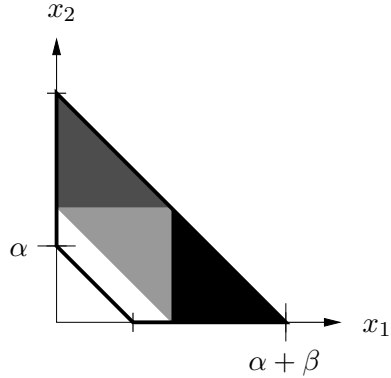}
 \end{center}
 \caption{Maximal packings of $(S^2 \ltimes S^2, \oo_{ab})$ by $2$
  and $3$ balls, $\frac{\bb}{\aa} \ge 1$.}  \label{figure48}
\end{figure}

\begin{figure}[h] 
 \begin{center}
  \psfrag{x1}{$x_1$}
  \psfrag{x2}{$x_2$}
  \psfrag{aa}{$\aa$}
  \psfrag{ab}{$\aa + \bb$}
  \psfrag{a}{$(a)$}
  \psfrag{b}{$(b)$}
  \psfrag{c}{$(c)$}
  \psfrag{12}{$\frac{\bb}{\aa} = \frac{1}{2}$}
  \psfrag{11}{$\frac{\bb}{\aa} = 1$}
  \psfrag{32}{$\frac{\bb}{\aa} = \frac{3}{2}$}
  \leavevmode\epsfbox{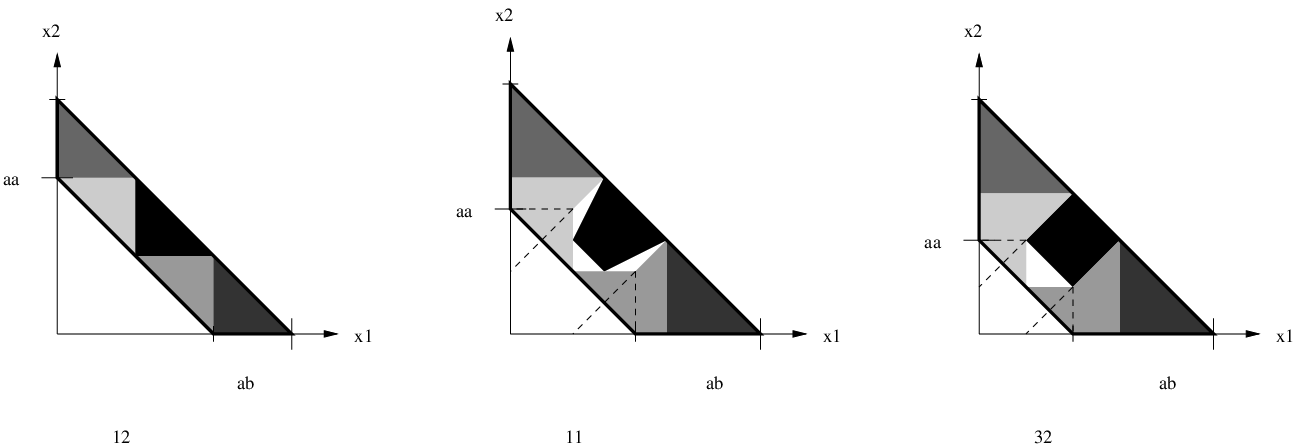}
 \end{center}
 \caption{Maximal packings of $(S^2 \ltimes S^2, \oo_{ab})$ by $4$
          and $5$ balls, $\frac{1}{2} \le \frac{\bb}{\aa} \le
          \frac{3}{2}$.}            \label{figure49}
\end{figure}

\begin{figure}[h] 
 \begin{center}
  \psfrag{x1}{$x_1$}
  \psfrag{x2}{$x_2$}
  \psfrag{aa}{$\aa$}
  \psfrag{ab}{$\aa + \bb$}
  \psfrag{a}{$(a)$}
  \psfrag{b}{$(b)$}
  \psfrag{c}{$(c)$} 
  \psfrag{33}{$(a): \frac{\bb}{\aa} =3$}
  \psfrag{22}{$(b1): \frac{\bb}{\aa} =2$}
  \psfrag{66}{$(b2): \frac{\bb}{\aa} =6$}
  \leavevmode\epsfbox{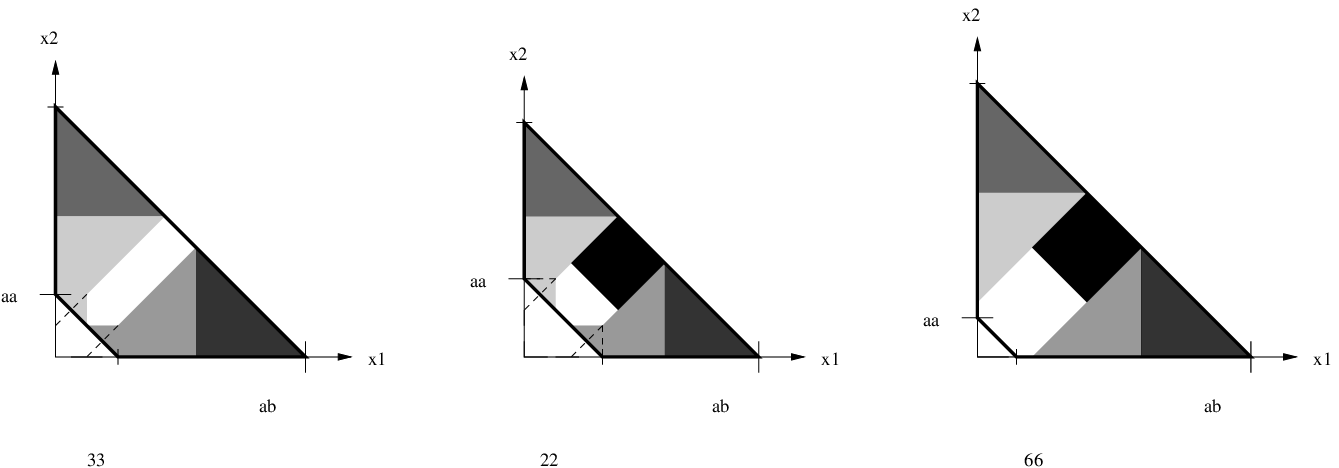}
 \end{center}
 \caption{Maximal packings of $(S^2 \ltimes S^2, \oo_{ab})$ by $4$
  and $5$ balls, $\frac{\bb}{\aa} \ge \frac{3}{2}$.}  \label{figure70}
\end{figure}

\smallskip
To construct packings with $p_k = k \frac{\bb}{2\aa+\bb}$ for all $k$ with 
$\lfloor \frac{k}{2} \rfloor \frac{\bb}{\aa} \le 1$ is very easy. 
\mbox{Figure \ref{figure47}(a)} shows a maximal packing by $1$ ball,
and \mbox{Figures \ref{figure47}(b1)} and (b2) show maximal packings by
$4$ and $5$ balls for $\frac{\bb}{\aa} = \frac{1}{2}$ and
$\frac{\bb}{\aa} < \frac{1}{2}$, respectively.  
%
%
\mbox{Figure \ref{figure48}} shows maximal packings for $k = 2,3$ and
$\frac{\bb}{\aa} \ge 1$. 
%
%
Also our maximal packings by $4$ balls are easy to understand
(Figure~\ref{figure49} and Figure~\ref{figure70}(a)):
$2 \, c_4 = \bb + \frac{\aa}{2}$ 
just means that
the two middle gray balls touch each other. As long as $\frac{\bb}{\aa} \le
\frac{3}{2}$, there is enough room for a fifth (black) ball between these two
balls. If $\frac{\bb}{\aa} > \frac{3}{2}$, there is enough space for a
fifth ball if and only if the capacity $c$ of the balls satisfies 
$2c + \frac{c}{2} \le \aa + \bb$; hence $c_5 = \frac{2\aa +2\bb}{5}$
(\mbox{Figures \ref{figure70}(b1)} and (b2)). 
%
%
%
%
%

\begin{figure}[h] 
 \begin{center}
  \psfrag{x1}{$x_1$}
  \psfrag{x2}{$x_2$}
  \psfrag{aa}{$\aa$}
  \psfrag{ab}{$\aa + \bb$}
  \psfrag{12}{$\frac{\bb}{\aa} = \frac{1}{2}$}
  \psfrag{23}{$\frac{\bb}{\aa} = \frac{2}{3}$}
  \psfrag{11}{$\frac{\bb}{\aa} = 1$}
  \leavevmode\epsfbox{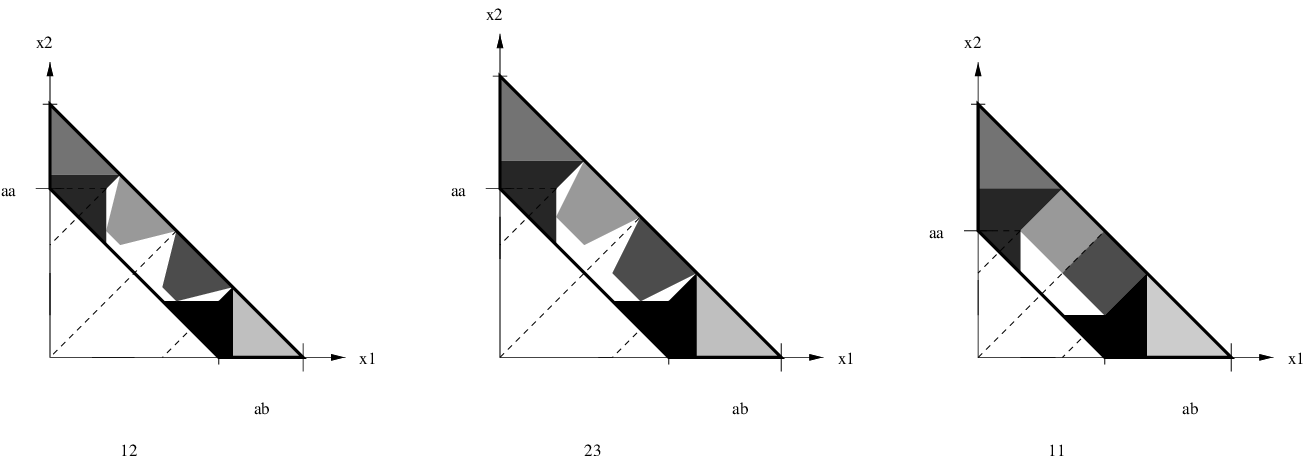}
 \end{center}
 \caption{Maximal packings of $(S^2 \ltimes S^2, \oo_{ab})$ by $6$
          balls, $\frac{1}{3}\le \frac{\bb}{\aa} \le 1$.}  \label{figure71}
\end{figure}

\begin{figure}[h] 
 \begin{center}
  \psfrag{x1}{$x_1$}
  \psfrag{x2}{$x_2$}
  \psfrag{aa}{$\aa$}
  \psfrag{ab}{$\aa + \bb$}
  \psfrag{a}{$(a)$}
  \psfrag{b}{$(b)$}
  \psfrag{c}{$(c)$} 
  \psfrag{44}{$\frac{\bb}{\aa} =4$}
  \psfrag{19}{$\frac{\bb}{\aa} =19$}
  \leavevmode\epsfbox{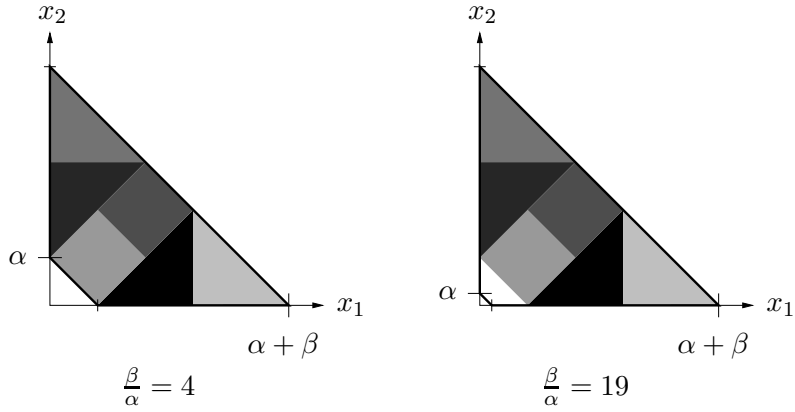}
 \end{center}
 \caption{Maximal packings of $(S^2 \ltimes S^2, \oo_{ab})$ by $6$
          balls, $\frac{\bb}{\aa} \ge 4$.}  \label{figure72}
\end{figure}

Let now $k=6$. Figure~\ref{figure71} shows maximal
packings for $\frac{1}{3} \le \frac{\bb}{\aa} \le 1$. For
$\frac{\bb}{\aa} > \frac{1}{3}$ the vertices of the ``lower middle
ball'' are
\[
(\aa + \bb - 2c_6, c_6), 
\quad \left( \frac{\aa + \bb}{2}, \frac{\aa + \bb}{2} \right), 
\quad (\aa + \bb -c_6, c_6), 
\quad \left( \frac{\aa + \bb}{2}, \frac{\aa + \bb}{2}-c_6 \right).
\]
%
%
Maximal packings for $\frac{\bb}{\aa} \ge 4$ are illustrated in
Figure~\ref{figure72}. 
%
%

\begin{remark}  \label{r:gap}
{\rm
It is not a coincidence that we were not able to construct
maximal packings of $\left( S^2 \ltimes S^2, \oo_{ab} \right)$
by $6$ balls for all ratios $\frac \bb \aa >0$.
Indeed, a maximal packing of $\left( S^2 \ltimes S^2, \oo_{ab}
\right)$ by $6$ equal balls for $\frac \bb \aa = \frac 53$ corresponds
to a maximal packing of the $4$-ball by $7$ equal balls.
\diam
}
\end{remark}

\begin{figure}[h] 
 \begin{center}
  \psfrag{x1}{$x_1$}
  \psfrag{x2}{$x_2$}
  \psfrag{aa}{$\aa$}
  \psfrag{ab}{$\aa + \bb$}
  \leavevmode\epsfbox{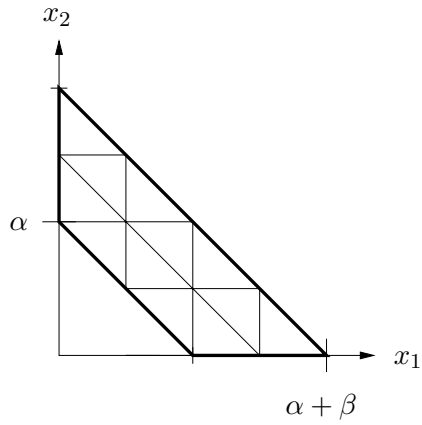}
 \end{center}
 \caption{A full packing of $(S^2 \ltimes S^2, \oo_{ab})$,
  $\frac{\bb}{\aa} =1$, by $12$ balls.}  \label{figure73}
\end{figure}

Finally, suppose that $\frac{\bb}{\aa} = \frac{l}{m-l}$ for some $l,m \in
\NN$ with $m>l$. We can then fill $(S^2 \ltimes S^2, \oo_{ab})$
by $l(2m-l)$ balls by decomposing $S^2 \ltimes S^2$ into $l$ shells
and filling the $i$-th shell with $2m+1-2i$ balls 
(see Figure~\ref{figure73}, where $l=2$ and $m=4$). 

%
%

\subsubsection{Maximal packings of $\Sigma_g(a) \times S^2(b)$ and 
$\left( \Sigma_g \ltimes S^2, \omega_{ab} \right)$ for $g \ge 1$.}
\label{223}

Fix $a>0$ and $b>0$.
We represent the symplectic structure of $\Sigma_g(a) \times S^2(b)$
by a split form.
Removing a wedge of $2g$ loops from $\Sigma(a)$ and a point from
$S^2(b)$ we see that $\Sigma_g(a) \times S^2(b)$ contains 
$\sq (a,b) \times \sq^2(1)$.
The explicit construction of the ``standard K\"ahler form'' in class $[
\oo_{ab}]$ given in \cite[Section~3]{M1} 
and \cite[Exercise~6.14]{MS}
shows that also $\left( \Sigma_g \ltimes S^2, \omega_{ab} \right)$
endowed with this standard form contains $\sq (a,b) \times \sq^2(1)$.
The explicit maximal packings promised after Proposition~\ref{p:C24} can thus
be constructed as for $S^2(a) \times S^2(b)$, see Figures~\ref{figure43}
and \ref{figure46}.

\subsection{Explicit packings of $\Sigma_g(a) \times \Sigma_h(b)$ for
$g,h \ge 1$} 
\label{333}

We consider $4$-manifolds of the form $\Sigma_g \times \Sigma_h$
with $g,h \ge 1$.
The space of symplectic structures on such manifolds in not understood,
but no symplectic structure different from 
$\Sigma_g(a) \times \Sigma_h(b)$ for some $a>0$, $b>0$ is known.
For $\Sigma_g(a) \times \Sigma_h(b)$, no obstructions to full packings
are known.
Recall from \eqref{P:finite} that for $\frac ab \in \QQ$,
\[
P \left( \Sigma_g(a) \times \Sigma_h(b) \right) \,:=\, 
\inf \left\{ k_0 \in \NN \mid p_k \left( \Sigma_g(a) \times \Sigma_h(b)
\right) =1 \quad \mbox{for all 
}\, k \ge k_0 \right\} 
\]
is finite.
In fact, Biran showed in Corollary~1.B and Section~5 of \cite{B2} that
\begin{equation}  \label{est9:T2}
P \left( T^2(1) \times T^2(1) \right) \,\le\, 2
\end{equation}
and that  
\begin{equation}  \label{est9:8ab2ab}  
P \left( \Sigma_g(a) \times \Sigma_h(b) \right) \,\le\, 
\left\{
  \begin{array}{ccl}   
         8 ab  & \text{ if }\, &  a,b \in \NN ,  \\ [.0em]
         2 ab  & \text{ if }\, &  a,b \in \NN \setminus \{ 1 \}.
  \end{array}   
\right.                   
\end{equation}
If $\frac ab \notin \QQ$ or if $1 \le k < P \left( \Sigma_g(a) \times
\Sigma_h(b) \right)$, not much is known about $p_k \left( \Sigma_g(a) \times \Sigma_h(b) \right)$:
We can assume without loss of generality that $a \ge b$. Since the
symplectic packing numbers of $S^2(a) \times S^2(b)$ and 
$\sq (a,b) \times \sq^2(1)$ agree, and since $\sq (a,b) \times \sq^2(1)$
symplectically embeds into $\Sigma_g(a) \times \Sigma_h(b)$,
\begin{equation}  \label{e9:abab}
p_k \left( S^2(a) \times S^2(b) \right) \,\le\, 
p_k \left( \Sigma_g(a) \times \Sigma_h(b) \right) \quad\, \text{ for all }\;
k \in \NN ,
\end{equation}
and Figures~\ref{figure43}, \ref{figure44}, \ref{figure45} and
\ref{figure46} describe some explicit packings of 
$\Sigma_g(a) \times \Sigma_h(b)$.
A comparison of Corollary~\ref{c:C21} with the estimates~\eqref{est9:T2} and
\eqref{est9:8ab2ab} and with Proposition~\ref{p9:jiang} below
shows, however, that in general the inequalities~\eqref{e9:abab} are not
equalities and that for $\Sigma_g(a) \times \Sigma_h(b)$ not all of 
the packings in Figures~\ref{figure43},
\ref{figure44} and \ref{figure45} are maximal.

Elaborating an idea of Polterovich, \cite[Exercise~12.4]{MS}, 
Jiang constructed in \cite[Corollary 3.3 and 3.4]{J}
explicit symplectic embeddings of one ball which improve the estimate
$\frac{b}{2a} \le p_1 \left( \Sigma_g(a) \times \Sigma_h(b) \right)$
from \eqref{e9:abab}.

\begin{proposition} 
{\rm (Jiang)}  \label{p9:jiang}
Let $\SS (a)$ be any closed surface of area $a \ge 1$.
\begin{itemize}
\item[(i)]
There exists a constant $C>0$ such that 
$p_1 \left( \SS(a) \times T^2(1) \right) \ge C$.
\item[(ii)]
If $h \ge 2$, there exists a constant $C=C(h) >0$ depending only on $h$ such
that $w_G \left( \SS(a) \times \SS_h (1) \right) \ge C \log a$.
In other words, 
\[
p_1 \left( \SS(a) \times \SS_h(1) \right) \,\ge\,
\frac{\left( C \log a \right)^2}{2a} .
\]
\end{itemize}
\end{proposition}

Notice that for $\SS = S^2$ Biran's result 
$p_1 \left( S^2(a) \times \SS_h(1) \right) = \min \left( 1, \frac{a}{2}
\right)$ stated in Proposition~\ref{p:C24} is much stronger.
We shall use Jiang's embedding method to prove the following
quantitative version of Proposition~\ref{p9:jiang}\:(i).

\begin{proposition}  \label{p9:jt}
If $a \ge 1$,
\[
p_1 \left( \SS(a) \times T^2(1) \right) \,\ge\, 
\frac{\max \{ a+1 - \sqrt{2a+1},2 \}}{4a}.
\]
In particular, the constant $C$ in Proposition~\ref{p9:jiang}\:(i) can be
chosen to be $C=1/8$.
\end{proposition} 

\proof
Set $R(a) = \{ (x,y) \in \RR^2 \, | \, 0<x<1, \; 0<y<a \}$, and consider the
linear symplectic map
\begin{eqnarray*}
\ff \colon  (R(a) \times R(a), dx_1 \wedge dy_1 + dx_2 \wedge dy_2)  &  \ra   &
                (\RR^2 \times \RR^2, dx_1 \wedge dy_1 + dx_2 \wedge dy_2),  \\
      (x_1,y_1,x_2,y_2) & \mapsto & (x_1+y_2,y_1, -y_2, y_1+x_2).
\end{eqnarray*}
Let $\pr \colon \RR^2 \ra T^2 = \RR / \ZZ \times \RR
/ \ZZ$ be the projection onto the standard symplectic torus. Then
$\left( \id_2 \times \pr \right) \circ \ff \colon R(a) \times R(a) \ra 
 \RR^2 \times T^2$ is a symplectic embedding. 
Indeed, given $(x_1,y_1,x_2,y_2)$ and $(x_1',y_1',x_2',y_2')$ with 
\begin{eqnarray}
x_1 + y_2 & =      & x_1' + y_2'           \label{equationI}  \\
      y_1 & =      &        y_1'           \label{equationII} \\
     -y_2 & \equiv & -y_2'        \mod \ZZ \label{equationIII}\\
y_1 + x_2 & \equiv & y_1' + x_2'  \mod \ZZ \label{equationIV} 
\end{eqnarray}
equations~\eqref{equationII} and \eqref{equationIV}
imply $x_2 \equiv x_2' \mod \ZZ$, whence $x_2 = x_2'$. Moreover,
\eqref{equationIII} and \eqref{equationI} show that $y_2-y_2' = x_1'-x_1
\equiv 0 \mod \ZZ$, and so $x_1 = x_1'$ and $y_2 = y_2'$.
Next observe that 
\[
\bigl( \left( \id_2 \times \pr \right) \circ \ff \bigr) \bigl( R(a) \times
R(a) \bigr) 
\,\subset\,\;
]-a,0[ \times ]-a-1, a+1[ \times T^2. 
\]
Thus $R(a) \times R(a)$ symplectically embeds into $\SS(2a(a+1)) \times
T^2(1)$, and since $B^4(a)$ symplectically embeds into $R(a) \times R(a)$
and $B^4(1)$ symplectically embeds into $\SS(a) \times T^2(1)$ 
for any $a \ge 1$, Proposition~\ref{p9:jt} follows.
\proofend

\subsection{Maximal packings of $4$-dimensional ellipsoids} 
\label{s9:34}

We finally construct some explicit maximal packings of $4$-dimensional
ellipsoids 
\[
E(a,b) \,=\,  
\left\{ (z_1, z_2 ) \in \CC^2 \; \bigg| \,
\frac{\pi |z_1|^2}{a} + \frac{\pi |z_2|^2}{b} < 1 \right\}.
\]
Without loss of generality we consider $E(\pi,a)$ with $a \ge \pi$.
\begin{proposition}  \label{p9:ellipsoid}
(i)
For each $k \in \NN$ the ellipsoid $E(\pi, k\pi)$ admits an
explicit full symplectic packing by $k$ balls.

\s
(ii)
$p_1 \left( E(\pi,a) \right) = \frac{\pi}{a}$ and 
$p_2\left( E(\pi,a) \right) = \min \left( \frac{2\pi}{a}, \frac{a}{2\pi}
\right)$, and these packing numbers can be realized by explicit
symplectic packings.
\end{proposition}
The statement~(i) was proved in \cite[Theorem~6.3\:(2)]{T}, 
and (ii) was proved in 
\cite[Corolary~3.11]{MMT}.
Their embeddings are different from ours.

\b
\ni
{\it Proof of Proposition~\ref{p9:ellipsoid}:} \,
(i) 
Set
\begin{eqnarray*}
\tr (a,b) &=& \left\{ 0 < x_1, x_2 \, \,\Big|\, 
\frac{x_1}{a}+\frac{x_2}{b} <1 \right\} \,\subset\, \RR^2(x), \\ [.2em]
\sq (a,b) &=& \{ 0 < y_1 < a, \,0<y_2<b  \} \,\subset\, \RR^2(y), 
\end{eqnarray*}
and abbreviate $\tr^2(\pi) = \tr(\pi,\pi)$.
We see as in Section~\ref{s:construction.1} that we can think of $B^4(\pi)$ as
$\tr^2(\pi) \times \sq^2(1)$ and of $E(\pi,k\pi)$ as $\tr (\pi,k\pi)
\times \sq^2(1)$.
The linear symplectic map 
$\left( x_1, x_2, y_1, y_2 \right) \mapsto \left( x_1, kx_2, y_1, \frac
1k y_2 \right)$ maps $\tr^2 (\pi) \times \sq^2(1)$ to
$\tr \left(\pi, k\pi\right) \times \sq \left( 1, \frac 1k \right)$, and
it is clear how to insert $k$ copies of this set into $\tr \left(
\pi,k\pi \right) \times \sq^2(1)$.

\s
(ii)
The estimates $p_1 \left( E(\pi,a) \right) \le \frac{\pi}{a}$ and 
$p_2 \left( E(\pi,a) \right) \le \frac{2\pi}{a}$ follow from the
inclusion $E(\pi,a) \subset Z^4(\pi)$ and from Gromov's Nonsqueezing
Theorem, 
and
$p_2\left( E(\pi,a) \right) \le \frac{a}{2\pi}$ follows from 
$E(\pi,a) \subset B^4(a)$ and Gromov's result 
$p_2 \left( B^4(a) \right) \le \frac 12$ stated in \eqref{eq9:pkBkn}.
The inclusion $B^4(\pi) \subset E(\pi,a)$ shows that $p_1 \left(
E(\pi,a) \right) = \frac \pi a$,
and explicit symplectic packings of $E(\pi,a)$ by two balls realizing
$p_2 \left( E(\pi,a) \right) = \min \left( \frac{2\pi}{a},
\frac{a}{2\pi} \right)$ can be constructed as in the proof of (i).
\proofend

\section{Maximal packings in higher dimensions}  \label{s9:higher}

\ni
In dimensions $2n \ge 6$, only few maximal symplectic packings by equal
balls are known.

\m
\ni
{\bf 1. Balls and $\left( \CC \PP^n, \oo_{SF} \right)$}

\s
\ni
As in dimension $4$ we denote by $\oo_{SF}$ the unique 
$\mbox{U}(n+1)$-invariant K\"ahler form on 
$\CC \PP^n$ whose integral over $\CC \PP^1$ equals $1$. 
The embedding \eqref{emb:BCP2} generalizes to all dimensions, and
\[
p_k \left( B^{2n} \right) \,=\, p_k \left( \CC \PP^n,
\oo_{SF} \right)
\quad\; \text{for all }\, k ,
\]
see \cite[Remark~2.1.E]{MP}.
Recall from \eqref{eq9:pkBkn} and \eqref{eq9:pkBln} that
\begin{eqnarray*}     
   p_k \left( B^{2n} \right) &=& \frac{k}{2^n}  \quad\,
   \text{ for }\, 2 \le k \le 2^n,   \\ [.3em]
   p_{l^n} \left( B^{2n} \right) &=& 1  
        \quad\, \text{ for all }\, l \in \NN .
\end{eqnarray*}
An explicit maximal packing of $\left( \CC \PP^n, \oo_{SF} \right)$ by 
$k \le n+1$ balls was found by Karshon in \cite{Ka}, 
and explicit full packings of $B^{2n}$ by $l^n$ balls for each $l
\in \NN$ were given by Traynor in \cite{T}.
Taking $l=2$, any $k$ balls of such a packing yield a maximal packing by
$k$ balls. The following different construction of an explicit full
packing of $B^{2n}$ by $l^n$ equal balls is mentioned in 
\cite[Remark~5.13]{T}.
Set
\begin{eqnarray*}
\tr^n (a) &=& \left\{ 0 < x_1, \dots, x_n \, \,\Bigg|\, \sum_{i=1}^n
\frac{x_i}{a} <1 \right\} \,\subset\, \RR^n(x), \\ [.3em]
\sq^n (a) &=& \{ 0 < y_i < a, \; 1 \le i \le n \}  \,\subset\, \RR^n(y) .
\end{eqnarray*}
We see as in Section~\ref{s:construction.1} that we can think of $B^{2n}(\pi)$
as $\tr^n (\pi) \times \sq^n(1)$ and of 
$B^{2n} \left( \frac \pi l \right)$ as $\tr^n \left( \frac \pi l \right)
\times \sq^n(1)$.
The matrix $\diag \left[ l, \dots, l, \frac 1l, \dots, \frac 1l \right]
\in \sym$ maps 
$\tr^n \left( \frac \pi l \right) \times \sq^n(1)$ to 
$\tr^n \left( \pi \right) \times \sq^n \left( \frac 1l \right)$.
It is clear how to insert $l^n$ copies of 
$\tr^n \left( \pi \right) \times \sq^n \left( \frac 1l \right)$ into 
$\tr^n \left( \pi \right) \times \sq^n \left( 1 \right)$.

\m
\ni
{\bf 2. Products of balls, complex projective spaces and surfaces}

\s
\ni
Set $n = \sum_{i=1}^d n_i$ and let $a_1, \dots, a_d \in \pi\NN$.
According to \cite[Theorem~1.5.A]{MP}, the product 
\[
\left( \CC \PP^{n_1} \times \dots \times \CC \PP^{n_d}, 
a_1\;\!\oo_{SF} \oplus \dots \oplus a_d \;\!\oo_{SF} \right)
\]
admits a full symplectic packing by
$\frac{n !}{n_1! \,\cdots\, n_d!} a_1^{n_1} \cdots a_d^{n_d}$ equal
$2n$-dimensional balls.
These full packings can be constructed in an explicit way. Indeed, 
explicit full packings of $B^{2n_i}(a_i)$ by $a_i^{n_i}$ equal balls as in
1.\ above can be used to construct explicit full packings of
\[
B^{2n_1}(a_1) \times \dots \times B^{2n_d}(a_d)
\]
by $\frac{n !}{n_1! \,\cdots\, n_d!} a_1^{n_1} \cdots a_d^{n_d}$ balls, see 
\cite[Section~3.2]{Kr}.
In particular, there are explicit full packings of the polydisc $P
\left( a_1, \dots, a_n \right)$ and of the products of surfaces 
$\Sigma_{g_1} (a_1) \times \cdots \times \Sigma_{g_n} (a_n)$ with 
$a_i \in \pi\NN$
by $n! \,a_1 \cdots a_n$ equal balls, see also \cite[Section~4.1]{T},
\cite[Theorem~4.1]{MMT}, and Figure~\ref{figure46} above for the case $n=2$.
An explicit packing construction in \cite[Theorem~1.21]{MMT} yields 
the lower bounds 
\[
p_7 \left( C^6 (\pi) \right) \ge \frac{224}{375}
\qquad \text{ and } \qquad
p_8 \left( C^6 (\pi) \right) \ge \frac{9}{16} .
\]
The technique in the proof of Proposition~\ref{p9:jt} can be used to
generalize Proposition~\ref{p9:jiang}\:(i): 
For any closed surface $\Sigma$ endowed with an area form $\ss$ and
any constant symplectic form $\oo$ on the $2n$-dimensional torus $T^{2n}$, 
there exists a constant $C>0$ such that
$p_1 \left( \Sigma \times T^{2n}, a \ss \oplus \oo \right) \ge C$ for
all $a \ge 1$, see \cite[Theorem~3.1]{J}.

\m
\ni
{\bf 3. Ellipsoids}

\s
\ni
Generalizing Proposition~\ref{p9:ellipsoid}\:(ii), the packing
numbers  
$p_1 \left( E(a_1, \dots, a_n) \right) = \frac{a_1^n}{a_1 \cdots a_n}$
and
$p_2 \left( E(a_1, \dots, a_n) \right) = \frac{2}{a_1 \cdots a_n} \min
\left( a_1^n, \left( \frac{a_n}{2} \right)^n\right)$
of a $2n$-dimensional ellipsoid were computed and realized by explicit
symplectic packings in 
\cite[Corollary~3.11]{MMT}.

\begin{remark}  \label{r9:explicit6}
{\rm
Karshon's explicit packing of $\left( \CC \PP^n, \oo_{SF} \right)$ by $k
\le n+1$ balls is maximal in the sense of \eqref{e9:full}.
Since in dimensions $\ge 6$ it is not yet known whether the space of
symplectic embeddings of a closed ball into a larger ball is connected,
all other explicit (and non-explicit) maximal symplectic packings known in
dimensions $\ge 6$ are maximal only in the sense of \eqref{e9:fullee}.
\diam
}
\end{remark}

\m
We conclude with addressing two widely open problems. As
before, we consider connected symplectic manifolds of finite volume.
\begin{question}  \label{q9:pk}
Which connected symplectic manifolds $(M,\oo)$ of finite volume satisfy
$p_k(M,\oo)=1$ for all $k \ge 1$?
\end{question}
Examples are $2$-dimensional manifolds, 
$\left( \CC \PP^2, \oo_{SF} \right)$ symplectically blown up at $N \ge 9$
points with weights close enough to $1/ \sqrt{N}$
and $S^2(1) \times S^2(1)$ symplectically blown up at $N \ge 8$
points with weights close enough to $1/ \sqrt{N}$ (see
\cite[Section~5]{B1}),
the ruled symplectic $4$-manifolds $\Sigma_g(a) \times S^2(b)$ and $\left(
\Sigma_g \ltimes S^2, \oo_{ab} \right)$ with $g \ge 1$ and $b \ge 2a$
and their symplectic blow-ups (see Proposition~\ref{p:C24} and
\cite[Theorem~6.A]{B1}), 
as well as certain closed symplectic $4$-manifolds
described in \cite[Theorem~2.F]{B1} and their symplectic
blow-ups.

A related problem is
\begin{question}  \label{q9:p1}
Which connected symplectic manifolds $(M,\oo)$ of finite volume satisfy
$p_1(M,\oo)=1$?
\end{question}
Examples different from the above ones are the ball $B^{2n}$ and
$\left( \CC \PP^n, \oo_{SF} \right)$, and, more generally, the
complement $\left( \CC \PP^n \setminus \Gamma, \oo_{SF} \right)$ of a
closed complex submanifold $\Gamma$ of $\CC \PP^n$ 
(see \cite[Corollary~1.5.B]{MP}).

\b
\ni
{\bf Acknowledgements.}
These pictures were painted in 1998 at ETH Z\"urich.
The words were added in 2004 at Leipzig University.
I thank both institutions for their support and  
Peter Albers for valuable discussions.


\markboth{{\rm References}}{{\rm References}}

\enddocument